\documentclass[12pt]{amsart}

\usepackage{amsmath, amsthm}
\usepackage{amssymb}
\usepackage{amscd}
\usepackage{latexsym}
\usepackage{epsfig}
\usepackage{graphics}
\usepackage{amsfonts}
\usepackage{psfrag}
\usepackage{graphicx}

\input xy
\xyoption{all}
\UseComputerModernTips

\numberwithin{equation}{section}
\numberwithin{figure}{section}

\newtheorem{theorem}{Theorem}[section]
\newtheorem{lemma}[theorem]{Lemma}
\newtheorem{proposition}[theorem]{Proposition}
\newtheorem{corollary}[theorem]{Corollary}

\theoremstyle{definition}

\theoremstyle{remark}

\newtheorem{remark}[theorem]{Remark}

\newcommand{\C}{{\mathbb{C}}}

\newcommand{\Z}{{\mathbb{Z}}}
\newcommand{\Q}{{\mathbb{Q}}}

\newcommand{\bbS}{\mathbb{S}}

\newcommand{\cupprod}{\smallsmile}

\newcommand{\hcup}{\smallsmile}

\newcommand{\kcup}{\cdot}
\newcommand{\kcap}{\smallfrown}
\newcommand{\less}{\!\setminus\!}
\newcommand{\dirac}{\mbox{$\not\negthinspace\partial$}}

\renewcommand{\t}{\mathfrak{t}}
\newcommand{\g}{\mathfrak{g}}
\newcommand{\h}{\mathfrak{h}}

\newcommand{\into}{\hookrightarrow}
\newcommand{\onto}{\twoheadrightarrow}

\renewcommand{\mod}{/\!/}

\newcommand{\pit}{p}

\DeclareMathOperator{\Tor}{Tor}
\DeclareMathOperator{\sign}{sgn}
\DeclareMathOperator{\pt}{pt}
\DeclareMathOperator{\Thom}{Thom}
\DeclareMathOperator{\Spinc}{Spin^{c}}
\DeclareMathOperator{\Index}{Index}
\DeclareMathOperator{\rank}{rank}
\DeclareMathOperator{\vertical}{Vert}
\DeclareMathOperator{\im}{im}
\DeclareMathOperator{\Ann}{ann}
\DeclareMathOperator{\Ad}{Ad}
\DeclareMathOperator{\Spin}{Spin}
\DeclareMathOperator{\SO}{SO}
\DeclareMathOperator{\U}{U}

\begin{document}

\title{The $K$-theory of abelian versus nonabelian symplectic quotients}

\author{Megumi Harada}
\address{Department of Mathematics, Bahen Centre for Information Technology, 40 St.\ George St., Room 6290, University of Toronto, Toronto, ON M5S2E4, Canada}
\email{megumi@math.toronto.edu} 
\urladdr{http://www.math.toronto.edu/\~{}megumi/}

\author{Gregory D. Landweber}
\address{Mathematics Department, University of Oregon,
Eugene, Oregon 97403-1222, USA}
\email{greg@math.uoregon.edu}
\urladdr{http://math.uoregon.edu/\~{}greg/}


\dedicatory{Dedicated to Raoul Bott}
\date{\today}



\begin{abstract} 
We compare the $K$-theories of symplectic quotients with respect to a
compact connected Lie group and with respect to its maximal torus,
and in particular we give a method for computing the former in terms
of the latter. More specifically, let $G$ be a compact connected Lie
group with no torsion in its fundamental group, let $T$ be a maximal
torus of $G$, and let $M$ be a compact Hamiltonian $G$-space. Let \(M
\mod G\) and \(M \mod T\) denote the symplectic quotient of $M$ by $G$
and by $T$, respectively. Using Hodgkin's K\"unneth spectral sequence
for equivariant $K$-theory, we express the $K$-theory of $M\mod G$ in
terms of the elements in the $K$-theory of $M\mod T$ which are
invariant under the action of the Weyl group, in addition to the Euler class
$e$ of a natural $\Spinc$ vector bundle over $M\mod T$.
This Euler class $e$ is induced by the denominator in the Weyl
character formula, viewed as a virtual representation of $T$; this 
is relevant for our proof.

Our results are $K$-theoretic analogues of similar (unpublished)
results by Martin for rational cohomology. However, our results and approach differ
from his in three significant ways. First, Martin's method involves
integral formul\ae, but the corresponding index formul\ae\ in
$K$-theory are too coarse a tool, as they cannot detect
torsion. Instead, we carefully analyze related $K$-theoretic
pushforward maps. Second, Martin's method involves dividing by the
order of the Weyl group, which is not possible in (integral) $K$-theory. We
render this unnecessary by examining Weyl anti-invariant elements,
proving a $K$-theoretic version of a lemma due to Brion. Finally,
Martin's results are expressed in terms of the annihilator ideal of
$e^{2}$, the square of the Euler class mentioned above. We are able
to ``remove the square,'' working instead with the annihilator ideal of
$e$.

\end{abstract}

\maketitle

\tableofcontents

\section{Introduction}

The main result of this paper is a formula for the $K$-theory of a
symplectic quotient $M \mod G$ with respect to a compact connected
Lie group $G$, in terms of the $K$-theory of the symplectic quotient $M\mod T$ with respect to a
maximal torus $T$ of $G$.
Such a formula is useful because computations of topological
invariants, such as cohomology or $K$-theory, in equivariant
symplectic geometry are frequently easier to carry out for {\em abelian}
group actions. Thus, it is a common practice in equivariant symplectic
geometry to prove results involving nonabelian compact connected group
actions by first reducing to simpler computations for their maximal
tori. Our theorem is an example of this general strategy. 

Symplectic quotients arise naturally in various fields. For example, toric varieties (which arise in combinatorics) and moduli spaces of bundles over Riemann surfaces (studied in gauge theory) are symplectic quotients. Furthermore, symplectic quotients are closely related to K\"ahler quotients
or Geometric Invariant Theory (``GIT'') quotients in complex algebraic geometry.
Hence, many moduli spaces that arise as GIT quotients also have symplectic realizations, and the topological invariants of such moduli spaces give useful constraints on moduli problems. In addition, the theory of geometric quantization
provides
a fundamental link between the topology of symplectic quotients and representation theory (see e.g. \cite[Section 7]{Kir98}). 



How, then, does one compute the cohomology or $K$-theory of symplectic quotients? 
In the case of rational cohomology, a fundamental result in this direction is the Kirwan surjectivity theorem \cite{Kir84}. 
Briefly, this theorem states that there is a ``model ring'' which
surjects onto the ring of interest; namely, there is a natural surjective ring
homomorphism 
\begin{equation}\label{eq:Kirwan-original}
\kappa: H^*_G(M;\Q) \onto H^*(M \mod G; \Q),
\end{equation}
where $H^*_G(M;\Q)$ is the $G$-equivariant
cohomology ring of the original Hamiltonian $G$-space from which
$M\mod G$ is constructed. 
In order to compute
$H^*(M\mod G;\Q)$, it therefore suffices to compute two objects: the
equivariant cohomology ring $H^*_G(M;\Q)$ and the kernel of
$\kappa$. The crucial advantage of this strategy is that for both of
these latter computations, the presence of the $G$-action on $M$
allows us to use \textit{equivariant} techniques which are
unavailable on the quotient.
Subsequent work 
(see e.g. \cite{JK95, GKM, TW03, Gol02}) 
has yielded a rich
theory for computing the {\em rational} cohomology rings of symplectic
quotients. However, satisfactory techniques are not known for {\em
integral} cohomology. 
We argue in \cite{HL05} that one
should, instead, look to the $K$-theory 
of
symplectic quotients as the natural integral setting in which to
extend the known results for rational cohomology. In particular, we
demonstrate that the extra torsion information 
in $K$-theory, unlike the torsion in integral cohomology, is amenable to
several of the standard techniques used in this field. These techniques
include the Atiyah-Bott lemma \cite{AB82}, which is the key tool for many
Morse-theoretic arguments in equivariant symplectic geometry.  In \cite{HL05},
we prove the $K$-theoretic analogue of the Kirwan surjectivity theorem
and examine the structure of $K_{G}^{*}(M)$, thereby
opening the door to computations of the $K$-theory of symplectic
quotients using techniques similar to those above.


Our main result in this manuscript is a $K$-theoretic analogue of an
(unpublished) result of Martin \cite{Mar00} for rational cohomology. This rational cohomology result also appears with an alternative proof due to Jeffrey, Mare, and Woolf in \cite{JMW02}.
We streamline Martin's original proof for rational cohomology in a manner which extends to $K$-theory, and in particular allows for the presence of torsion. 
Thus, by using our more direct methods, we are able to go beyond
computations of Betti numbers to subtler invariants of the
manifolds. We will discuss in more detail, below, the differences
between our approach and Martin's approach. We also note that in \cite{HP05}, Hausel and Proudfoot have developed a hyperk\"ahler version of Martin's theorem, 
but for this paper we restrict our attention to the symplectic case.

Martin's result is not a
direct computation of the kernel of the Kirwan map, but
is instead a comparison between a kernel of a nonabelian Kirwan
map and an abelian Kirwan map. We now recall the
setting of Martin's, and therefore our, results in some more
detail. Throughout this paper, let $G$ be a compact connected Lie
group, let $T$ be a maximal torus in $G$, and let $W = N(T)/T$ denote
the corresponding Weyl group. For our purposes, we impose only one
torsion constraint, that $G$ have no torsion in its fundamental group
$\pi_{1}(G)$. On the other hand, if there is torsion in $\pi_{1}(G)$, then we have two possible recourses: first, we could eliminate the torsion by replacing $G$ by a finite cover, or second, we could work with $K$-theory with coefficients in a ring $S$ for which $\pi_{1}(G)\otimes S$ is torsion free, i.e., inverting those primes which appear in the torsion. We note that Martin requires no such assumptions about the torsion of
$G$ since he works with rational or real cohomology.

Given a compact Hamiltonian $G$-space $M$ with moment map $\mu_G: M
\to \g^*$ such that $0$ is a regular value of $\mu_G$, the symplectic
quotient $M \mod G$ is defined to be the quotient \(\mu_G^{-1}(0)/G.\)
Here, the quotient is viewed as the standard quotient if $G$ acts
freely on $\mu_G^{-1}(0)$, or an orbifold otherwise. Given $T$ a
maximal torus in $G$, the inclusion \(\t \into \g\) of the
corresponding Lie algebras gives rise to the linear projection \(\g^* \to
\t^*,\) and composing this with $\mu_G$ yields the $T$-moment map
$\mu_T: M \to \t^*$, with \(M \mod T := \mu_T^{-1}(0)/T.\) Here again
we assume \(0 \in \t^*\) is a regular value of $\mu_T$, thus making
$M\mod T$ a manifold if the $T$-action on $\mu_T^{-1}(0)$ is free, and an orbifold
otherwise. 
Below, we assume for simplicity that both $M\mod G$ and $M\mod T$
are manifolds, but our arguments readily generalize to the orbifold case.

In this setting, the inclusions \(\mu_G^{-1}(0) \into M\) and \(\mu_T^{-1}(0) \into M\) induce natural maps
\begin{equation}\label{eq:Kirwan-maps}
\kappa_G: K^*_G(M) \to K^*(M \mod G), \quad \kappa_T: K^*_T(M) \to K^*(M \mod T), 
\end{equation}
which we show in \cite{HL05} are surjective.
The problem of comparing $K^*(M \mod G)$ to $K^*(M \mod T)$ can therefore be formulated in terms of comparing the kernels of $\kappa_G$ and $\kappa_T$. 
We have the following theorem:

\begin{theorem}\label{theorem:martin-intro}
	Let $G$ be a compact connected Lie group with no torsion in
	its fundamental group $\pi_{1}(G)$, and let $T$ be a maximal
	torus in $G$. If $M$ is a compact Hamiltonian $G$-space, and $0$
	is a regular value of the moment maps $\mu_{G}$ and $\mu_{T}$, then
	the kernels of the Kirwan maps $\kappa_{G}$ and $\kappa_{T}$ given by (\ref{eq:Kirwan-maps}) are related by
	$$\ker \kappa_{G}
	\cong \left\{ x\in K^{*}_{T}(M)^{W} : \kappa_{T}(x)\kcup e = 0 \right\},$$
where $e \in K^{*}(M\mod T)$ is the class induced by the denominator of the Weyl character formula under the map
$$R(T) \cong K_{T}^{*}(\pt) \to K_{T}^{*}\bigl(\mu_{T}^{-1}(0)\bigr) \cong K^{*}(M\mod T).$$
\end{theorem}

We first make some observations about this class $e$ which appears in
the statement of the above theorem.  The denominator of the Weyl
character formula is given by
\begin{equation*}
	\prod_{\alpha > 0} \left( e^{\alpha/2} - e^{-\alpha/2} \right)
	= \sum_{w\in W} \sign(w) \, e^{w(\rho)},
\end{equation*}
in terms of a choice of positive root system for $G$, where $\rho$ is half the sum of the positive roots. The choice of positive root system affects only the overall sign of the Weyl denominator, but since we are interested in the annihilator of the induced element $e\in K^{*}(M\mod T)$, this overall sign does not make any difference to Theorem~\ref{theorem:martin-intro}. Furthermore, even if we were to factor out an overall $e^{\rho}$ from the Weyl denominator, as is done in some treatments of the Weyl character formula, this still does not affect the annihilator appearing in Theorem~\ref{theorem:martin-intro}.

We now observe that the Weyl denominator is in fact an equivariant
$K$-theoretic Euler class, and hence arises naturally in our geometric
context.  If $\g$ and $\t$ are the Lie algebras of $G$ and $T$,
respectively, then the Weyl denominator is in fact the $T$-equivariant
Euler class
$e_{T}(\g/\t) \in R(T)$ of the representation $\g/\t$, viewed as a
$T$-equivariant bundle over a point. More geometrically, the Weyl
denominator can also be viewed as the $G$-equivariant Euler class
$e_{G}(G/T) \in K_{G}^{*}(G/T)$ of the tangent bundle to the
homogeneous space $G/T$, under
the isomorphism $K_{G}^{*}(G/T) \cong R(T)$. The overall
sign and possible factor of $e^{\rho}$ in the choice of Euler class
corresponds to the choice of $\Spinc$-structure required to construct
a Thom isomorphism in $K$-theory, and is analogous to the choice of
orientation in cohomology.


Our second main theorem states our result in terms of the $K$-theories
of the two symplectic quotients. This second version is slightly
stronger, since its proof requires the surjectivity of the Kirwan map
$\kappa_G$ in~\eqref{eq:Kirwan-maps} above.

\begin{theorem}
Let $G$ be a compact connected Lie group with no torsion in its fundamental group $\pi_{1}(G)$, and let $T$ be a maximal torus in $G$. If $M$ is a compact Hamiltonian $G$-space, and $0$
	is a regular value of the moment maps $\mu_{G}$ and $\mu_{T}$, then the $K$-theories of the symplectic quotients \(M \mod G\) and \(M \mod T\) are related by an isomorphism
\begin{equation}
K^*(M\mod G) \cong \frac{K^*(M \mod T)^W}{\Ann(e)},
\end{equation}
where $e\in K^{*}(M\mod T)$ is as described in Theorem~\ref{theorem:martin-intro}.
\end{theorem}

Our theorems and methods differ in three significant ways from
that of the
rational or de~Rham cohomology versions given in \cite{JMW02} and
\cite{Mar00}. First, the proofs of the cohomological versions rely on
integral formul\ae\ which do not work in $K$-theory in the presence of
torsion.
%
 In particular, integral formul\ae\ involve real-valued pairings coming from Poincar\'e duality, and the corresponding index formul\ae\ in $K$-theory involve integer-valued pairings, which necessarily vanish on elements of finite order. In this paper, instead of using integral formul\ae\, we carefully analyze the $K$-theory pushforward map $i_{!}:K^{*}(X) \to K^{*}(Y)$ induced by a particular $\Spinc$-inclusion $i : X \hookrightarrow Y$. Secondly, the proofs in \cite{JMW02} and \cite{Mar00} involve dividing by $|W|$, the order of the Weyl group, either in the statements of the integral formul\ae\ themselves, or in the process of averaging over the Weyl group. In an integral theory such as $K$-theory, making sense of such formul\ae\ or techniques requires inverting those primes which divide $|W|$.
In place of these arguments, we consider Weyl anti-invariance, proving the following $K$-theoretic version of a lemma which Brion uses in \cite{Bri91} to establish a related result for the cohomology of geometric invariant theory quotients:
\begin{lemma}
	If $X$ is a compact $G$-space, then the Weyl anti-invariant elements $K_{T}^{*}(X)^{-W}$ are a free module over $K_{G}^{*}(X) \cong K_{T}^{*}(X)^{W}$ with generator $e$ as described in Theorem~\ref{theorem:martin-intro}.
\end{lemma}

This lemma appears as Part (2) of our Corollary~\ref{cor:brion} below, where we show that it is a consequence of the Weyl character formula, as interpreted by  Bott in \cite{Bott65}. Finally, the statement of our theorems differs from Martin's \cite[Theorem A]{Mar00}, since he uses the annihilator ideal of $e^{2}$ in place of the annihilator ideal of the Euler class $e$.
The distinction between $e$ and $e^{2}$ and the equivalence of the two corresponding statements is proved and discussed in detail in \cite{JMW02}.

We now briefly survey the contents of this paper. In
Section~\ref{sec:background}, we recall some essential
symplecto-geometric and topological facts relating nonabelian and
abelian symplectic quotients. Then in Section~\ref{sec:weyl}, we
develop the key tools which allow us to carefully analyze properties
of Weyl invariance and anti-invariance in equivariant $K$-theory. We
define a $K$-theoretic pushforward for a $\Spinc$-inclusion in
Section~\ref{sec:pushforward}, and also prove a key lemma involving
multiplication by the Euler class $e$. We prove our two main theorems
in Section~\ref{sec:Ktheory-MmodG-MmodT}. The last two sections
discuss generalizations of our results, first to the case where we
replace $T$ by an arbitrary subgroup $H$ of maximal rank in $G$, and
secondly to other generalized cohomology theories.

Finally, since the premise of this manuscript is that it is ``easier''
to compute the kernel of the Kirwan map $\kappa_T$ for abelian $T$, we
conclude the introduction with a few comments regarding the
computation in the abelian case. 
The arguments used, for example, in \cite{TW03} to give
an explicit description of the kernel of $\kappa_T$ in the case of
rational cohomology depend essentially
on the (rational-cohomological) Atiyah-Bott lemma. 
By
using the $K$-theoretic version of the Atiyah-Bott lemma \cite{HL05,
VeVi03} and the results of this manuscript, it is straightforward
to generalize to the $K$-theoretic setting the explicit description of the kernel of the Kirwan map
(in rational cohomology) given by Tolman and Weitsman \cite{TW03}.
We also expect some generalization of
Goldin's effective algorithm \cite{Gol02} for computing the kernel in the abelian
case to also hold in $K$-theory. These topics, together with examples computing the $K$-theory of symplectic quotients, will be the subject of a future paper \cite{HL-kernel}.

\bigskip

\noindent {\bf Acknowledgements.}  It is our pleasure to thank Lisa
Jeffrey, Peter Landweber, Nick Proudfoot, Dev Sinha, Peter Teichner, and Jonathan Weitsman for helpful discussions, and Victor Guillemin and Allen Knutson for their interest and encouragement. The second author thanks the University of
Toronto and the Fields Institute for their hospitality and support
while conducting a portion of this research. Both authors thank the
Banff International Research Station for their hospitality.

\section{Background}\label{sec:background}
 

Let $G$ be a compact, connected, possibly nonabelian Lie group. Suppose that
$M$ is a compact Hamiltonian $G$-space. Let $T$ denote a
choice of maximal torus of $G$; then $M$ is also naturally a
Hamiltonian $T$-space, by restricting the action of $G$ and composing the $G$-moment map with the natural projection of Lie coalgebras.
As mentioned in the introduction, our main goal in this paper is to relate the
$K$-theory of the nonabelian symplectic quotient $M \mod G$ to that of the abelian symplectic quotient $M \mod T$. In this section, we recall the basic facts which make our computation possible. 

We begin by reviewing the topological relationship between the two
symplectic quotients. 
Denote by $\mu_G$ the
moment map for the $G$-action. Then the corresponding moment map $\mu_T$ for the $T$-action
is related to $\mu_G$ by the diagram
\[
\xymatrix{
K^{*}_{G}(M) \ar[dr]_{\mu_T} \ar[r]^-{\mu_G} & \g^* \ar[d]^{\pit} \\
 & \t^*, \\
}
\]
where $\pit : \g^* \to \t^*$ is the linear projection dual to the inclusion \(\t \into \g.\)
Throughout this paper, we assume that $0$ is a regular value of both $\mu_{G}$ and $\mu_{T}$. As a consequence, the groups $G$ and $T$ act locally freely on the level sets $\mu_{G}^{-1}(0)$ and $\mu_{T}^{-1}(0)$, respectively. The corresponding symplectic quotients or Marsden-Weinstein reductions of $M$ at $0$ are then
$$M\mod G := \mu_{G}^{-1}(0) / G, \qquad
  M\mod T := \mu_{T}^{-1}(0) / T.$$
We view these as standard quotients if the group actions are free, or as orbifolds otherwise. 
In the orbifold case, our $K$-theoretic statements should be interpreted in terms of the orbifold $K$-theory: 
$$K_{\text{orb}}^{*}(M\mod G) := K^{*}_{G}\left(\mu_{G}^{-1}(0)\right),
\qquad
K_{\text{orb}}^{*}(M\mod T) := K^{*}_{T}\left(\mu_{T}^{-1}(0)\right).$$
In other words, in the orbifold case, 
we work directly with the equivariant $K$-theory of the level sets instead of dealing with the ordinary $K$-theory of the quotient. 
Our arguments in Section~\ref{sec:weyl} are already equivariant and are therefore unaffected by this difference, and the arguments in Section~\ref{sec:pushforward} have straightforward generalizations to the equivariant case. Thus, henceforth, we assume for simplicity that $M \mod G$ and $M \mod T$ are manifolds.



The definition of the symplectic quotient in both the nonabelian and
the abelian case immediately gives rise to two natural maps, induced
by the equivariant inclusions of the zero sets of the moment maps,
\(\mu_G^{-1}(0) \into M\) and \(\mu_T^{-1}(0) \into M:\) 
\begin{equation}\label{eq:Kirwan-nonabelian}
\xymatrix{
K_{G}^{*}(M) \ar[r] \ar[dr]_{\kappa_{G}} & K_{G}^{*}\left(\mu_{G}^{-1}(0)\right) \ar[d]^{\cong} \\
& K^{*}(M\mod G)
} \qquad
\xymatrix{
K_{T}^{*}(M) \ar[r] \ar[dr]_{\kappa_{T}}& K_{T}^{*}\left(\mu_{T}^{-1}(0)\right) \ar[d]^{\cong} \\
& K^{*}(M\mod T)
}
\end{equation} 
Both $\kappa_G$, $\kappa_T$ are called {\bf Kirwan maps.} 
In \cite{HL05}, we prove that both $\kappa_G$ and $\kappa_T$ are
{\em surjective.} Thus, in order to compute either $K^*(M \mod G)$ or
$K^*(M \mod T)$, it suffices to compute the kernels of the respective
Kirwan maps. In other words, we have 
\begin{equation}
K^*(M \mod G) \cong \frac{K^*_G(M)}{\ker(\kappa_G)}, \quad K^*(M \mod
T) \cong \frac{K^*_T(M)}{\ker(\kappa_T)}. 
\end{equation}
Furthermore, there is a natural forgetful map 
\begin{equation}\label{eq:forgetful}
r^{G}_T: K^*_G(M) \to K^*_T(M),
\end{equation}
obtained by restricting the $G$-action on a bundle to a $T$-action, which Atiyah shows is split injective in \cite[Proposition 4.9]{Ati68}.
We have already mentioned in the introduction that there are
many known methods for explicitly computing the kernel $\ker(\kappa_T)$ for an
abelian group $T$. The main result of this paper may now be stated as
follows: given an element \(\alpha \in K^*_G(M),\) we will show that
\(\alpha \in \ker(\kappa_G)\) if and only if $r^{G}_T(\alpha) \in K^*_T(M)$ satisfies a
certain condition related to $\ker(\kappa_T)$. Thus, we will reduce
the question of computing the kernel $\ker(\kappa_G)$ to that of computing the
kernel $\ker(\kappa_T)$. 


 We first note that since \(\mu_T = \pit \circ \mu_G,\) and \(0 \in \g^*\) maps to \(0 \in \t^*\) under $\pit$, we have that \(\mu_G^{-1}(0) \subseteq \mu_T^{-1}(0).\) Moreover, the action of $T$ preserves both $\mu_G^{-1}(0)$ and $\mu_T^{-1}(0)$, so this is a $T$-equivariant inclusion. Hence the submanifold $\mu_G^{-1}(0)/T$ includes into the $T$-symplectic quotient \(M \mod T := \mu_T^{-1}(0)/T.\) On the other hand, the nonabelian symplectic quotient \(M \mod G\) is defined as \(\mu_G^{-1}(0)/G.\) Thus we have the following diagram relating the two symplectic quotients, using the intermediate space \(\mu_G^{-1}(0)/T:\) 
\begin{equation}\label{eq:topological-Martin}
\xymatrix{
	\mu_{G}^{-1}(0) / T \,\ar@{^(->}[r]^-{i} \ar@{->>}[d]_{\pi} & \mu_{T}^{-1}(0) / T = M\mod T \\
	\mu_{G}^{-1}(0) / G = M\mod G
}
\end{equation}
where the vertical arrow $\pi$ is a fibration with fiber \(G/T.\) 

Combining the Kirwan maps from (\ref{eq:Kirwan-nonabelian}), the forgetful map (\ref{eq:forgetful}), and the maps induced by (\ref{eq:topological-Martin}), we obtain the following commutative diagram:
\begin{equation}
\label{eq:Kirwan-GtoT-diagram}
\xymatrix{
K^{*}_{G}(M) \ar@{->>}[d]_{\kappa_{G}} \, \ar@{^(->}[rr]^{r^{G}_{T}} & & K^{*}_{T}(M) \ar@{->>}[d]^{\kappa_{T}} \\
K^{*}(M\mod G) \, \ar@{^(->}[r]^-{\pi^{*}} & K^{*}\left(\mu_{G}^{-1}(0)/T\right) & K^{*}(M\mod T)\, .\ar[l]_-{i^{*}}
}
\end{equation}
Here we recall that the Kirwan maps $\kappa_{G}$ and $\kappa_{T}$ are surjective by our prior result \cite{HL05}, and that the forgetful map $r^{G}_{T}$ is split injective by \cite{Ati68}. In addition, the map $\pi^{*}$ can be identified with the forgetful map $r^{G}_{T}:K^{*}_{G}(\mu^{-1}_{G}(0)) \to K^{*}_{T}(\mu^{-1}_{G}(0))$, and thus 
is injective as well. By a diagram chase, we see that
\begin{equation}\label{eq:kernels}
\ker \kappa_{G} = \ker (\pi^{*}\circ \kappa_{G}) = \ker (i^{*}\circ \kappa_{T}\circ r^{G}_{T}) \cong \ker (i^{*}\circ \kappa_{T})|_{\text{im}\,r^{G}_{T}}.
\end{equation}

We now recall the following key topological fact \cite[Proposition 1.2]{Mar00} about the diagram~\eqref{eq:topological-Martin} which will allow us to describe the kernel of $i^*$ explicitly.
Let $\t^{0}\subset \g^{*}$ denote the annihilator of $\t$ in $\g^{*}$, or equivalently let $\t^{0} = \ker \pit$ be the kernel of the projection $\pit : \g^{*}\to \t^{*}$.
Consider the moment map $\mu_{G} : M \to \g^{*}$ restricted to  $\mu_{T}^{-1}(0)$.
Since $\mu_{T} = \pit \circ \mu_{G}$, we see that $\mu_{G}$ restricts to a $T$-equivariant map
\begin{equation}\label{eq:moment-map-t0}
\mu_{G} : \mu_{T}^{-1}(0) \to \t^{0}.
\end{equation}
Furthermore, $\mu_{G}^{-1}(0)$ is precisely the zero set of this restricted moment map. Viewing $\mu_{T}^{-1}(0)$ as a principal $T$-bundle over the symplectic quotient $M\mod T$, let
\begin{equation}\label{eq:bundle}
E := \mu_{T}^{-1}(0) \times_{T} \t^{0} \to M\mod T
\end{equation}
be the associated vector bundle  induced by the representation $\t^{0}$ of $T$. The restricted moment map (\ref{eq:moment-map-t0}) then gives a section $s : M\mod T \to E$ whose zero set is $\mu_{G}^{-1}(0)/T$.

\begin{proposition}\label{prop:topological-martin}
The vector bundle $E \to M\mod T$ given by (\ref{eq:bundle}) has a section $s$ transverse to the zero section whose zero set is $\mu_{G}^{-1}(0) / T$. Furthermore, if $G$ has no $2$-torsion in its fundamental group $\pi_{1}(G)$, then a choice of $\Ad$-invariant inner product on $\g$ induces a unique $T$-equivariant spin structure $\Spin(E) \cong \mu_{T}^{-1}(0) \times_{T} \Spin(\t^{0})$ on $E$.
\end{proposition}

\begin{proof}
For reference, we briefly recall the argument given by Martin \cite{Mar00}. To obtain the desired $s$, we consider $\mu_{G}$ restricted to $\mu_{T}^{-1}(0)$, as in (\ref{eq:moment-map-t0}). This restriction if $T$-invariant, so by taking a quotient, we obtain a section $s$ of \(E = \mu_T^{-1}(0) \times_T \t^0.\) The fact that $0$ is a regular value of $\mu_G$ implies that $0 \in \t^0$ is a regular value of $\mu_G |_{\mu_T^{-1}(0)}$. Thus $s$ is transverse to the zero section. 

To prove the second statement, note that an inner product on $\g$ allows us to identify the annihilator $\t^{0}$ with the quotient $\g/\t$. If the inner product is $\Ad$-invariant, then the action of $T$
on $\t^{0}$ gives us a group homomorphism $T \to \SO(\t^{0}) \cong \SO(\g/\t)$,
and a $T$-equivariant spin structure corresponds to a lift to a group homomorphism $T \to \Spin(\t^{0}) \cong \Spin(\g/\t)$. However, this is precisely the same information required for a $G$-equivariant spin structure on the homogeneous space $G/T$, and in particular, if $\pi_{1}(G)$ has no $2$-torsion, then $G/T$ admits a unique homogeneous spin structure (see \cite{Bott65}).
%
\end{proof}

\section{Equivariant $K$-theory and the Weyl group}\label{sec:weyl}

The following lemma is a result of McLeod \cite{McL79}, which we  use in Section~\ref{sec:Ktheory-MmodG-MmodT} to compare the equivariant $K$-theories with respect to a compact Lie group and its maximal torus.
We prove it here in detail due to its importance to this paper. See also our discussion in \cite{HL05}, in which we use this result in order to establish a version of equivariant formality for Hamiltonian $G$-spaces. We note that a related result for algebraic $K$-theory appears in \cite{Mer98}.
 
\begin{lemma}\label{lemma:mcleod}
	Let $G$ be a compact connected Lie group with no torsion in its fundamental group $\pi_{1}(G)$, and let $X$ be a compact $G$-manifold. If $T$ is a maximal torus of $G$, then $X$ is likewise a $T$-manifold by restriction, and its $G$- and $T$-equivariant $K$-theories are related by
	\begin{equation}\label{eq:mcleod}
	K^{*}_{T}(X) \cong K^{*}_{G}(X) \otimes_{R(G)} R(T).
	\end{equation}
\end{lemma}

\begin{proof}
This result is a consequence of Hodgkin's K\"unneth spectral sequence for equivariant $K$-theory \cite{Hod75}, which computes $K^{*}_{G}(X\times Y)$ for two $G$-spaces $X$ and $Y$ via an Eilenberg-Moore spectral sequence starting with $E_{2}$ page
\begin{equation}\label{eq:e2}
E_{2} = \Tor^{R(G)}_{*}\bigl( K^{*}_{G}(X), K^{*}_{G}(Y) \bigr).
\end{equation}
In particular, the $0$-torsion is simply the tensor product,
$$\Tor_{0}^{R(G)} \bigl( K^{*}_{G}(X), K^{*}_{G}(Y) \bigr) \cong
K^{*}_{G}(X) \otimes_{R(G)}K^{*}_{G}(Y).$$
 In \cite[\S 3]{Sna72}, Snaith argues that under a technical hypothesis, later verified by McLeod in \cite{McL79}, this spectral sequence always converges as expected, provided that $G$ is a compact connected Lie group with $\pi_{1}(G)$ torsion-free.

In our case, we take $Y = G/T$ and use the isomorphisms
$$K^{*}_{G}(G/T) \cong R(T), \qquad K^{*}_{G}(X \times G/T) \cong K^{*}_{T}(X).$$
The former is given in \cite{Seg68} and is induced by the map identifying a representation $U$ of $T$ with the associated homogeneous bundle $G \times_{T} U$, or conversely identifying a homogeneous bundle $E \to G/T$ with its fiber over the identity coset $E_{eT}$. The latter is a generalization:
at the level of
vector bundles, a $T$-equivariant bundle $\pi: E \to X$ induces a
$G$-equivariant bundle $G\times_{T}E$ over $G/T\times X$. Here $G$ acts on the left on the $G$-factor of $G \times_{T}E$, while it acts via the diagonal action on $G/T \times X$. The projection map is given by $(g,v) \mapsto \left(gT,\,g\cdot \pi(v)\right)$
for $g\in G$ and $v\in E$, which is well defined since for any $t\in T$,
$$\left(gt^{-1},\,tv\right) \sim (g,v) \longmapsto \left(gt^{-1}T,\,gt^{-1}t\cdot \pi(v)\right) = \left(gT,\,g\cdot \pi(v)\right),$$
and is $G$-equivariant since for any $h\in G$,
$$h\cdot(g,v) = (hg,v) \longmapsto \left(hgT,\,hg\cdot\pi(v)\right) = h\cdot\left(g,\,g\cdot\pi(v)\right).$$
Conversely, given a $G$-equivariant bundle over $G/T \times X$, its restriction
to the identity coset $eT \times M$ is a $T$-equivariant bundle over $M$. These maps of vector bundles are inverses of each other and extend to equivariant $K$-theory.

Finally, Pittie shows in \cite{Pit72} that $R(T)$ is a free $R(G)$-module, and so we have
$$\Tor_{p}^{R(G)} \bigl( K_{G}^{*}(X), R(T) \bigr) = 0 \text{ for $p\neq 0$.}$$
The K\"unneth spectral sequence therefore degenerates into a standard K\"unneth isomorphism,
$$K^{*}_{T}(X) \cong K^{*}_{G}(X \times G/T) \cong K^{*}_{G}(X) \otimes_{R(G)}K^{*}_{T}(G/T)
\cong K^{*}_{G}(X) \otimes_{R(G)}R(T),$$
which gives us our desired result.
\end{proof}

\begin{remark}

Taking $X = G$ in Lemma~\ref{lemma:mcleod}, we obtain
\begin{equation}\label{eq:gmodt-tensor}
K^{*}(G/T) \cong K^{*}_{T}(G) \cong K^{*}_{G}(G) \otimes_{R(G)} R(T) \cong \Z \otimes_{R(G)}R(T),
\end{equation}
where $R(G)$ acts on $\Z$ by its augmentation homomorphism taking a virtual $G$-module to its integral dimension. Since $R(T)$ is a free $R(G)$-module by \cite{Pit72}, we can write it in terms of generators $a_{1}, \ldots, a_{m}$, where $m = |W|$, as
\begin{equation}\label{eq:rt-module}
R(T) \cong R(G) \, a_{1} \oplus \cdots \oplus R(G) \, a_{m}.
\end{equation}
Considering the composition
$$\alpha : R(T) \to K^{*}_{G}(G/T) \to K^{*}(G/T)$$
together with the isomorphism (\ref{eq:gmodt-tensor}),
we see that $K^{*}(G/T)$ is the free $\Z$-module
\begin{equation}\label{eq:gmodt-module}
K^{*}(G/T) \cong \Z\,\alpha(a_{1}) \oplus \cdots \oplus \Z\,\alpha(a_{m}),
\end{equation}
whose generators over $\Z$ are the images of the generators of $R(T)$ over $R(G)$.
Using (\ref{eq:rt-module}) and  (\ref{eq:gmodt-module}), we can rewrite the isomorphism (\ref{eq:mcleod}), replacing the tensor product over $R(G)$ with a tensor product over $\Z$, giving us the following non-canonical isomorphism:
\begin{equation*}\begin{split}
K_{T}^{*}(X) &\cong
K_{G}^{*}(X) \otimes_{R(G)} R(T) \cong \bigoplus_{i=1}^{m} K_{G}^{*}(X) \otimes_{R(G)} R(G)\, a_{i} \\
&\cong \bigoplus_{i=1}^{m} K_{G}^{*}(X) \otimes_{\Z} \Z\, \alpha(a_{i})
\cong K_{G}^{*}(X) \otimes_{\Z} K^{*}(G/T)
\end{split}\end{equation*}
as modules over $K^{*}_{G}(X)$.

The corresponding statement in cohomology over a coefficient ring $S$ is the isomorphism
\begin{equation}\label{eq:factor-cohomology}
H^*_{T}(X;S) \cong H^*_{G}(X;S) \otimes H^{*}(G/T;S),
\end{equation}
as modules over $H^{*}_{G}(X;S)$. This can be derived from the Leray-Serre spectral sequence for the fibration $G/T \to X_{T} \to X_{G}$,
which, in the words of \cite{AB82}, behaves like a product in cohomology. However, in order for (\ref{eq:factor-cohomology}) to hold, the cohomology $H^{*}(G;S)$ must be torsion-free. In other words, all of the primes dividing the torsion in the integral cohomology $H^{*}(G;\Z)$ must be invertible in $S$. This torsion constraint in cohomology is stronger than in $K$-theory, where we require only that the fundamental group $\pi_{1}(G)$ be torsion-free.  We can similarly bypass this $K$-theory torsion constraint by considering $K$-theory with coefficients in a ring $S$ where $\pi_{1}(G)\otimes S$ is torsion-free, i.e., by inverting all primes which occur in the torsion of the fundamental group $\pi_{1}(G)$.

We note that (\ref{eq:factor-cohomology}) is the cohomological form of our Lemma~\ref{lemma:mcleod} which Brion uses in his related discussion \cite{Bri91} of geometric invariant theory quotients.
\end{remark}

We now recall two fundamental facts about the representation ring of a compact Lie group. Let $G$ be a compact, connected, possibly nonabelian Lie group, and let $T$ be a maximal torus.  The Weyl group $W = N(T)/T$ acts on representations of $T$, and the restriction to $T$ of any representation of $G$ is invariant under the action of the Weyl group. In fact, at the level of representation rings, the Weyl group acts on $R(T)$, and the restriction map
$r^{G}_{T}:R(G) \to R(T)$ gives an isomorphism
$$r^{G}_{T} : R(G) \xrightarrow{\cong} R(T)^{W},$$
identifying $R(G)$ with the Weyl-invariants in $R(T)$. This is a consequence of the classification of the irreducible representations of $G$ in terms of their highest weights, together with the fact that the Weyl group acts simply transitively on the Weyl chambers.

Secondly, the Weyl group is generated by reflections, and viewing it as a group of permutations of the Weyl chambers, we have a homomorphism $\sign : W \to \{\pm 1\}$. The Weyl anti-invariant part $R(T)^{-W}$ of $R(T)$ consists of all elements on which $w\in W$ acts by $\sign(w)$. In \cite{Bott65}, Bott points out that the Weyl character formula,
\begin{equation}\label{eq:character}
\chi ( V_{\lambda} ) = \frac{ \sum_{w\in W} \sign(w)\, e^{w(\lambda + \rho)} } {\sum_{w\in W} \sign(w)\, e^{w(\rho)}},
\end{equation}
implies that the Weyl anti-invariant part $R(T)^{-W}$, which is the span of the principal alternating elements appearing in the numerator of (\ref{eq:character}), is a free module over $R(G)$ generated by the denominator of (\ref{eq:character}). This Weyl denominator can also be written as the $T$-equivariant Euler class of the representation $\g/\t$, as we derive in (\ref{eq:e-gmodt}) below,
\begin{equation}\label{eq:denominator}
e_{T}(\g/\t) = \prod_{\alpha>0}\left(e^{\alpha/2} - e^{-\alpha/2}\right),
\end{equation}
in terms of a choice of positive roots of $G$. We can now use Lemma~\ref{lemma:mcleod} to extend these facts to equivariant $K$-theory by a simple argument: 

\begin{corollary}\label{cor:brion}
	Let $G$ be a compact connected Lie group with no torsion in its fundamental group $\pi_{1}(G)$, and let $X$ be a compact $G$-space. Let $T$ be a maximal torus and $W = N(T)/T$ the Weyl group.
\begin{enumerate}
	\item The restriction map $r^{G}_{T}: K^{*}_{G}(X) \to K^{*}_{T}(X)$ gives an isomorphism
	$$r^{G}_{T} : K_{G}^{*}(X) \xrightarrow{\cong} K_{T}^{*}(X)^{W}$$
	onto the Weyl-invariant elements.
	\item The Weyl anti-invariant elements $K^{*}_{T}(X)^{-W}$ are a free module over $K^{*}_{G}(X)$ via the restriction map $r^{G}_{T}$, with a single generator given by the image of $e_{T}(\g/\t)\in R(T)^{-W}$, the Weyl denominator (\ref{eq:denominator}), under the natural map $R(T)\to K^{*}_{T}(X)$.
\end{enumerate}
\end{corollary}

\begin{proof}
Taking the Weyl-invariant part of both sides of (\ref{eq:mcleod}), and noting that the Weyl group acts trivially on $K_{G}^{*}(X)$ and $R(G)$, we have
\begin{equation*}\begin{split}
K_{T}^{*}(X)^{W} &\cong K_{G}^{*}(X) \otimes_{R(G)} R(T)^{W}
\\ &\cong K_{G}^{*}(X) \otimes_{R(G)}R(G) \cong K_{G}^{*}(X),
\end{split}\end{equation*}
proving (1). Likewise, taking the Weyl anti-invariant part of both sides of (\ref{eq:mcleod}),
we have
\begin{equation*}\begin{split}
K_{T}^{*}(X)^{-W} &\cong K_{G}^{*}(X) \otimes_{R(G)} R(T)^{-W} \\
&\cong K_{G}^{*}(X) \otimes_{R(G)}R(G)\cdot e_{T}(\g/\t) \\ &\cong K_{G}^{*}(X) \cdot e_{T}(\g/\t),
\end{split}\end{equation*}
giving us (2).
\end{proof}

\begin{remark}
Part (1) of Corollary~\ref{cor:brion} appears in \cite{McL79}. On the
other hand, we have not encountered part (2) of Corollary~\ref
{cor:brion} elsewhere in the $K$-theory literature. Part (2) is the $K
$-theoretic version of a lemma due to Brion \cite{Bri91}, which
states that the Weyl anti-invariant elements of $H_{T}^{*}(X)$, with
rational coefficients, are a free module of rank $1$ over $H_{G}^{*}
(X)$, generated by the image of any non-zero element of the top
cohomology group $H^{\text{top}}(G/T)$ under the isomorphism (\ref{eq:factor-cohomology}).
\end{remark}

\section{$K$-theoretic pushforwards, pullbacks, and Euler classes}\label{sec:pushforward}

In this section, we develop further tools that we will need for our arguments in Section~\ref{sec:Ktheory-MmodG-MmodT}. The original argument given in \cite{Mar00} for the cohomological statement is based on integral formul\ae\  which we do not have at our disposal for $K$-theory due to the presence of torsion. Hence we must develop $K$-theoretic techniques in order to achieve the same goal. We will phrase the results in this section more generally than the specific case of \(\mu_G^{-1}(0)/T\) and \(\mu_T^{-1}(0)/T\) under consideration in this paper. Recall that \(\mu_G^{-1}(0)/T \into \mu_T^{-1}(0)/T\) is an inclusion whose normal bundle is equipped, by Proposition~\ref{prop:topological-martin}, with a spin structure, and hence a canonical $\Spinc$ structure. We now consider the general situation. 

Let $i : X \hookrightarrow Y$ be an inclusion of compact smooth manifolds, and suppose that the normal bundle $N \to X$ of this inclusion comes equipped with a $\Spinc$ structure. (In fact, all we really require is a $\Spinc$ structure on the stable normal bundle.)
Let
$$k : N \to Y, \qquad
  k : ( N, N\less X ) \to (Y, Y\less X)$$
be the identification of the normal bundle with a tubular neighborhood of $X$ in $Y$ and the corresponding map of pairs.
We define the pushforward in $K$-theory, as introduced in \cite
{AS68I}, to be the composition
$$ i_{!} : K^{*}(X) \xrightarrow{\Thom} K^{*}\bigl(N, \, N\less X \bigr) \xrightarrow{(k^{*})^{-1}} K^{*}\bigl( Y, \, Y \less X\bigr) \xrightarrow{j^{*}} K^{*}(Y)$$
of the Thom isomorphism for $N$, the inverse of an excision isomorphism, and the $K$-theory homomorphism induced by the map of pairs
$$j : Y \to ( Y, Y \less X ).$$
Explicitly, letting $\pi : N \to X$ be the vector bundle projection, the pushforward $i_{!}$ acts on elements $x\in K^{*}(X)$ by
\begin{equation}\label{eq:pushforward}
	i_{!} : x \mapsto j^{*} \circ (k^{*})^{-1}\bigl( \pi^{*}x \kcup u(N) \bigr),
\end{equation}
where
$$u(N) \in K^{*}\bigl(N,\,N \less X\bigr) \cong K^{*}\bigl(Y,\, Y \less X\bigr)$$
is the $K$-theory Thom class of $N$ determined by its $\Spinc$ structure.
The product in \eqref{eq:pushforward}
is the map on relative $K$-theory
$$ K^{*}(N) \otimes K^{*}\bigl(N,\,N\less X \bigr) \to K^{*}\bigl(N,\,N\less X \bigr)$$
induced by the tensor product of vector bundles. We note that this definition of the pushforward does not depend on the choice of identification of the normal bundle $N$ with a tubular neighborhood of $X$ in $Y$. However, it does depend on the $\Spinc$ structure on $N$, which enters through the Thom class $u(N)$.

The following lemma first appeared in \cite{BS58} in the context of the $K$-theory of coherent sheaves. Here we present a topological version.

\begin{lemma}\label{lemma:pushforward-module}
	The pushforward map $i_{!} : K^{*}(X) \to K^{*}(Y)$ in $K$-theory induced by a $\Spinc$-inclusion $i : X \hookrightarrow Y$ of smooth manifolds is a $K^{*}(Y)$-module homomorphism, i.e.,
	\begin{equation}\label{eq:pushforward-identity}
		i_{!}( i^{*}y \kcup x ) = y \kcup i_{!}x
	\end{equation}
	for all $x \in K^{*}(X)$ and $y\in K^{*}(Y)$.
\end{lemma}

\begin{proof}
The naturality of the product on relative $K$-theory gives us the commutative diagram:
$$
\xymatrix{
	K^{*}(Y) \otimes K^{*} \bigl( Y,\, Y\less X\bigr) \ar[r]^-{\kcup} \ar@<-3.5em>[d]^{k^{*}} \ar@<2em>[d]_{\cong}^{k^{*}} &
	K^{*}\bigl(Y , \, Y \less X \bigr) \ar[d]_{\cong}^{k^{*}}
	\\
	K^{*}(N) \otimes K^{*} \bigl( N,\, N \less X \bigr) \ar[r]^-{\kcup} &
	K^{*}\bigl( N, \, N \less X \bigr)
}$$
Since the composition $i\circ \pi : N \to X \hookrightarrow Y$ is homotopic to $k$, we have
$$ k^{*} (y \kcup z) = \pi^{*}i^{*}y \kcup k^{*}z, $$
for any  $y\in K^{*}(Y)$ and $z\in K^{*}(Y,Y\less X)$.
The pushforward of $i^{*}y \kcup x$ for $y\in K^{*}(Y)$ and $x\in K^{*}(X)$ is then given by
\begin{align*}
i_{!}(i^{*}y \kcup x)
&= j^{*} \circ (k^{*})^{-1} \Bigl( \pi^{*} ( i^{*}y \kcup x ) \kcup u(N) \Bigr) \\
&= j^{*} \circ (k^{*})^{-1} \Bigl( \pi^{*}i^{*}y \kcup \pi^{*} x \kcup u(N) \Bigr) \\
&= j^{*} \Bigl( \, y \kcup (k^{*})^{-1}\bigl( \pi^{*}x \kcup u(N) \bigr) \, \Bigr).
\end{align*}
On the other hand, the products for standard and relative $K$-theory are related by the commutative diagram
$$
\xymatrix{
	K^{*}(Y) \otimes K^{*} \bigl( Y,\, Y\less X\bigr) \ar[r]^-{\kcup} \ar@<-3em>@{=}[d] \ar@<2em>[d]^{j^{*}} &
	K^{*}\bigl(Y , \, Y \less X \bigr) \ar[d]^{j^{*}}
	\\
	K^{*}(Y) \otimes K^{*}(Y)
	\ar[r]^-{\kcup} &
	K^{*}(Y)
}$$
and so we obtain
$$ i_{!}(i^{*}y \kcup x) = y \kcup \Bigl( j^{*} \circ (k^{*})^{-1} \bigl( \pi^{*}x \kcup u(N) \bigr) \Bigr)
= y \kcup i_{!}x,$$
which gives us our desired identity \eqref{eq:pushforward-identity}.
\end{proof}

\begin{remark}
If both $X$ and $Y$ are $\Spinc$, then we can alternatively define this pushforward homotopically as the composition
$$ i_{!} : K^{*}(X) \xrightarrow{P.D.} K_{\dim X - *}(X) \xrightarrow{\,i_{*}\,}
K_{\dim X - *}(Y) \xrightarrow{P.D.} K^{* + \dim Y - \dim X}(Y)$$
of Poincar\'{e} duality on $X$, and the standard pushforward on $K$-homology, and Poincar\'{e} duality on $Y$.  Following \cite{Ati70}, we view the $K$-homology group $K_{0}(X)$ as being determined by classes of elliptic operators on $X$, where two such operators are equivalent if and only if their homogeneous symbol determines the same class in $K^{0}_{c}(T^{*}X)$ (using $K$-theory with compact supports). Recalling from \cite{AS68I} that the symbol map is surjective if we include elliptic pseudo-differential operators, we see that the symbol map gives an isomorphism
\begin{equation}\label{eq:k-homology-tx}
\sigma : K_{0}(X) \xrightarrow{\cong} K^{0}_{c}(T^{*}X).
\end{equation}
If $X$ is $\Spinc$, then the Poincar\'{e} duality map is the composition
$$P.D. : K^{*}(X) \xrightarrow{\Thom} K^{*+\dim X}(T^{*}X) \xrightarrow{\sigma^{-1}} K_{\dim X - *}(X)$$
of the Thom isomorphism for the cotangent bundle $T^{*}X$ with the extension of the isomorphism \eqref{eq:k-homology-tx} to all degrees of $K$-theory, and it corresponds geometrically to the map
$ x \mapsto [\dirac_{X}\otimes x] $
taking a $K$-theory class $x\in K(X)$ to the $K$-homology class of the $\Spinc$ Dirac operator on $X$ twisted by $x$, where $x$ is viewed as the class of a virtual bundle.
The pushforward $i_{*}$ in $K$-homology can be defined naturally in terms of $C^{*}$-algebras.

We note that the identity \eqref{eq:pushforward-identity} is a general property of the $K$-theory pushforward $i_{!}$. It is not limited to the case where $i$ is an inclusion and does not depend on our particular definition \eqref{eq:pushforward}. To see this, we need to consider the $K$-theory version of the cap product, which using the identification \eqref{eq:k-homology-tx} is determined by the product on the $K$-theory of $T^{*}X$ according to the commutative diagram: 
$$\xymatrix{
	K^{*}(X) \ar@<-5ex>[d]^{\pi^{*}}_{\cong} \otimes K_{*}(X) \ar@<5ex>[d]^{\sigma}_{\cong} \ar[r]^-{\kcap} & K_{*}(X) \ar[d]^{\sigma}_{\cong} \\
	K^{*}(T^{*}X) \otimes K_{c}^{*}(T^{*}X) \ar[r]^-{\kcup} & K_{c}^{*}(T^{*}X)
}$$
Geometrically, this $K$-theory cap product is induced by the map taking an elliptic operator $D$ and a vector bundle $E$ to the twisted operator $D \otimes E$ on $E$-valued sections, obtained by choosing a connection on $E$ and replacing the derivatives appearing in $D$ with covariant derivatives. In this notation, Poincar\'{e} duality is the map $K^{*}(X) \to K_{*}(X)$ given by $x \mapsto [\dirac_{X}] \kcap x$.
Applying Poincar\'{e} duality, Lemma~\ref{lemma:pushforward-module} is equivalent to the statement that the pushforward $i_{*} : K_{*}(X) \to K_{*}(Y)$ in $K$-homology is a $K^{*}(Y)$-module homomorphism. The identity \eqref{eq:pushforward-identity} is equivalent to the identity
$$ i_{*} ( i^{*} y \kcap x ) = y \kcap i_{*} x $$
for $y \in K^{*}(Y)$ and $x\in K_{*}(X)$ in terms of the $K$-theory cap product.
However, this is nothing more than the naturality of the cap product, as expressed by the following commutative diagram:
$$\xymatrix{
	K^{*}(X) \otimes K_{*}(X) \ar[r]^-{\kcap} \ar@<4ex>[d]^{i_{*}} & K_{*}(X) \ar[d]^{i_{*}}\\
	K^{*}(Y) \otimes K_{*}(Y) \ar[r]^-{\kcap} \ar@<4ex>[u]^{i^{*}} & K_{*}(Y)
}$$
The corresponding statement holds for the cap product in homology and cohomology (see, for example, \cite[\S 3.3]{Hat02}).
\end{remark}

Composing the restriction map $i^{*}$ with the pushfoward $i_{!}$,
we obtain
$$ i^{*} \circ i_{!} : x \mapsto x \kcup e(N)$$
for all $x\in K^{*}(X)$,
as the Euler class $e(N)\in K^{*}(X)$ of the normal bundle $N$ is the restriction of the Thom class $u(N)$ to $X$. The reverse composition $i_{!}\circ i^{*}$ is generally not so simple.  However, in the case that we are interested in, we have an extension of $N$ to a bundle over all of $Y$, which allows us to compute this composition directly.

\begin{lemma}\label{lemma:transverse}
Given a $\Spinc$ vector bundle $E \to Y$ and a section $s : Y \to E$ transverse to the zero section, let $X$ be the zero set of $s$, and $i : X \hookrightarrow Y$ the inclusion. Then,
\begin{enumerate}
\item \label{enum:1}
The normal bundle $N$ to $X$ is isomorphic to the restriction: $N \cong i^{*}E$.
\item \label{enum:2}
The unit element $1\in K^{*}(X)$ pushes forward to the Euler class $i_{!}1 = e(E) \in K^{*}(Y)$.
\item \label{enum:3}
For all $y\in K^{*}(Y)$, we have the identity $i_{!}\circ i^{*}: y \mapsto y \kcup e(E)$.
\end{enumerate}
\end{lemma}

\begin{proof}
Statement (\ref{enum:1}) follows immediately from the transversality of the section $s$ determining $X$. Statement (\ref{enum:3}) follows from statement (\ref{enum:2}), together with Lemma~\ref{lemma:pushforward-module}, by noting that
$$i_{!} ( i^{*} y ) = i_{!} ( i^{*} y \kcup 1 ) = y \kcup i_{!}1.$$
It remains to show Statement (\ref{enum:2}). In light of our definition \eqref{eq:pushforward} of the pushforward, we must show that the Thom class $u(N) \in K^{*}(N,N\less X)$ maps to the Euler class $e(E)\in K^{*}(Y)$ via the composition
\begin{equation}\label{eq:k-j-composition}
	K^{*}(N,N\less X) \xrightarrow{(k^{*})^{-1}}
   	K^{*}(Y,Y\less X) \xrightarrow{\phantom{()}j^{*}\phantom{{}^{-1}}} K^{*}(Y) .
\end{equation}
Consider the Thom class $u(E) \in K^{*}(E,E\less Y)$ of the bundle $E \to Y$. By Statement (\ref{enum:1}) and the naturality of the Thom class, we have
\begin{equation}\label{eq:u-for-N-iE}
	u(N) = u(i^{*}E) = i^{*} u(E).
\end{equation}
On the other hand, the homomorphism
\begin{equation}\label{eq:Thom-to-Euler}
	K^{*}(E, E\less Y) \xrightarrow{l^{*}} K^{*}(E) \cong K^{*}(Y),
\end{equation}
pulling back to $Y$ via any section, maps the Thom class to the Euler class,
$ u(E) \mapsto e(E) $.

Since $X$ is the zero set of the section $s : Y \to E$, we have a map of pairs $$s : ( Y, Y\less X ) \to ( E, E \less Y ).$$ Combining that with the inclusion $$i : ( i^{*}E, i^{*}E\less X ) \to ( E, E\less Y)$$  and the isomorphism $N \cong i^{*}E$ from Statement (\ref{enum:1}), we obtain the commutative diagram
$$
\xymatrix{
	K^{*}(i^{*}E, i^{*}E\less X) \ar[d]_{\cong} & \ar[l]_-{i^{*}} K^{*}( E, E \less Y ) \ar[r]^-{l^{*}} \ar[d]^{s^{*}}
	& K^{*}(E) \ar[d]^{s^{*}}_{\cong} \\
	K^{*}(N, N\less X) & \ar[l]_-{k^{*}}^{\cong} K^{*}( Y, Y \less X ) \ar[r]^-{j^{*}} & K^{*}(Y) 
}
$$
(We note that the maps $i : N \cong i^{*}E \to E $ and $s \circ k : N \to E$ 
are homotopic when restricted to a sufficiently small neighborhood of the zero section $X$ inside of $N$, and so they induce the same maps on the relative $K$-theory.)
Starting with the Thom class $u(E) \in K^{*}(E,E\less Y)$, we obtain the commutative diagram:
$$
\xymatrix{
	u(i^{*}E) \ar@{|->}[d]_{\cong} & \ar@{|->}[l]_-{i^{*}} u(E) \ar@{|->}[r]^-{l^{*}} \ar@{|->}[d]^{s^{*}}
	& l^{*}u(E) \ar@{|->}[d]^{s^{*}}_{\cong} \\
	u(N) & \ar@{|->}[l]_-{k^{*}}^-{\cong} s^{*}u(E) \ar@{|->}[r]^-{j^{*}} & e(E) 
}
$$
where the left and right sides are given by \eqref{eq:u-for-N-iE} and  \eqref{eq:Thom-to-Euler} respectively.
This diagram shows that the composition \eqref{eq:k-j-composition} on the bottom row indeed maps the Thom class $u(N)$ to the Euler class $e(E)$.
\end{proof}

\begin{corollary}\label{corollary-index}
	Under the conditions of Lemma~\ref{lemma:transverse}, and further assuming that $Y$ is $\Spinc$, we have the index formula:
	\begin{equation}\label{eq:index}
	\Index \dirac_{X} \otimes i^{*} y
	= \Index \dirac_{Y} \otimes \bigl( y \kcup e(E) \bigr)
	\end{equation}
	for all $y\in K^{*}(Y)$. Here $\dirac_{Y}$ and $\dirac_{X}$ denote the Dirac operators corresponding to the given $\Spinc$ structure on $Y$ and the induced $\Spinc$ structure on $X$ respectively, and $\dirac \otimes x$ denotes the Dirac operator twisted by the $K$-theory class $x$.
\end{corollary}

\begin{proof}
We note that the following diagram commutes:
$$
\xymatrix{
	K^{*}(X) \ar[rr]^{i_{!}} \ar[dr]_{\Index \dirac_{X}} & & K^{*}(Y) \ar[dl]^{\Index\dirac_{Y}} \\
	& \Z
}
$$
and the index formula \eqref{eq:index} follows from the identity $i_{!}\,i^{*}y = y \kcup e(E).$
\end{proof}

\begin{remark}
The de~Rham cohomology analogue of the index formula \eqref{eq:index} is the integral formula
\begin{equation}\label{eq:integral}
	\int_{X} i^{*} y = \int_{Y} y \cupprod e_{H}(E),
\end{equation}
for all $y\in H^{*}_{dR}(Y)$, where $e_{H}(E) \in H^{\rank E}_{dR}(Y)$ is the cohomology Euler class.  In this case we require only that $Y$ and the bundle $E\to Y$ be oriented, rather than $\Spinc$. Alternatively, we could express the integral formula \eqref{eq:integral} in terms of evaluation on the top homology class,
\begin{equation}\label{eq:integral-homology}
	\bigl\langle i^{*}y, \,[X] \bigr\rangle
	= \bigl\langle y \cupprod e_{H}(E),\, [Y] \bigr\rangle,
\end{equation}
which also holds for cohomology with rational or integral coefficients.

The pushforward formula $i_{!}\,i^{*}y = y \kcup e(E)$ for $K$-theory has as its cohomological analogue the formula $i_{!}\,i^{*}y = y \cupprod e_{H}(E)$, where 
here $i_{!} : H^{*}(X) \to H^{*+\rank E}(Y)$ is the pushforward in cohomology. An equivalent statement is that the cohomology Euler class $e_{H}(E)$ of an oriented bundle $E\to Y$ is Poincar\'{e} dual to the zero set $X$ of a generic section transverse to the zero section (see \cite{BT82}). For de~Rham and rational cohomology, the pushforward formula is equivalent to the formul\ae\ \eqref{eq:integral} and \eqref{eq:integral-homology} respectively. However, for $K$-theory and integral cohomology, the index formula \eqref{eq:index} and formula \eqref{eq:integral-homology} respectively are weaker than the corresponding pushforward formul\ae. This is because the index of the Dirac operator and pairing with the top homology class are homomorphisms to the integers, and must therefore vanish on torsion elements.
\end{remark}

\section{The $K$-theory of $M\mod G$ versus $M \mod T$}\label{sec:Ktheory-MmodG-MmodT} 

We are now in a position to prove our two main theorems. 
Since we know by Part (1) of Corollary~\ref{cor:brion} that the restriction $r_T^G$ maps onto the Weyl-invariant elements, we may restrict the commutative diagram~\eqref{eq:Kirwan-GtoT-diagram} to the Weyl-invariant components of all the $K$-theory groups. This gives us the following new commutative diagram: 
%
\begin{equation}
\label{eq:Kirwan-GtoT-diagram-Weyl}
\xymatrix{
K^{*}_{G}(M) \ar@{->>}[d]_{\kappa_{G}} \ar[rr]^{r^{G}_{T}}_{\cong} & & K^{*}_{T}(M)^{W} \ar[d]^{\kappa_{T}} \\
K^{*}(M\mod G) \ar[r]^-{\pi^{*}}_-{\cong} & K^{*}\left(\mu_{G}^{-1}(0)/T\right)^{W} & K^{*}(M\mod T)^{W} \ar[l]_-{i^{*}}
}
\end{equation}
%
%
We recall from \cite{HL05} that the Kirwan maps $\kappa_{G}$ and $\kappa_{T}$ given by (\ref{eq:Kirwan-nonabelian}) are surjective. Since the Kirwan map $\kappa_{T}$	is essentially the restriction to $\mu_{T}^{-1}(0)$, we observe that it commutes	with the action of the Weyl group, and restricting to the Weyl-invariants, we obtain
\begin{equation}\label{eq:kirwan-weyl}
\kappa_{T} : K_{T}^{*}(M)^{W} \to K^{*}(M\mod T)^{W}.
\end{equation}
Similarly, the pullback $i^{*}$ restricts to give a map
$$i^{*} : K^{*}(M\mod T)^{W} \to K^{*}\bigl(\mu_{G}^{-1}(0)/T\bigr)^W$$
on the Weyl-invariant elements.

\begin{remark}
The Kirwan map~\eqref{eq:kirwan-weyl} restricted to Weyl-invariants is not necessarily surjective, as there may be Weyl-invariant elements in $K^{*}(M\mod T)$ which are the image of elements of $K_{T}^{*}(M)$ which are not themselves Weyl-invariant. If we were to work over the rationals,
then we could average over the Weyl group to construct Weyl-invariant preimages, but that may not be possible over the integers. This issue does not affect either of the two theorems below.
\end{remark}

We can now prove the main result of this article, which we state in two different forms. First, we state it in terms of the kernels of the Kirwan maps:

\begin{theorem}\label{theorem:martin}
	Let $G$ be a compact connected Lie group with no torsion in its fundamental group $\pi_{1}(G)$, and let $T$ be a maximal torus in $G$. If $M$ is a compact Hamiltonian $G$-space, and $0$
	is a regular value of the moment maps $\mu_{G}$ and $\mu_{T}$, then the kernels of the Kirwan maps $\kappa_{G}$ and $\kappa_{T}$ given by (\ref{eq:Kirwan-nonabelian}) are related by
	$$\ker \kappa_{G}
	\cong \left\{ x\in K^{*}_{T}(M)^{W} : \kappa_{T}(x)\kcup e(E) = 0 \right\},$$
where $e(E)\in K^{*}(M\mod T)$ is the $K$-theoretic Euler class of the bundle $E$ given by (\ref{eq:bundle}).
\end{theorem}

Second, we state it in terms of the $K$-theories of the symplectic quotients. This second version is slightly stronger, and it assumes our Kirwan surjectivity result from \cite{HL05}.

\begin{theorem}\label{theorem:Ktheory-Martin} 
Let $G$ be a compact connected Lie group with no torsion in its fundamental group $\pi_1(G)$, and let $T$ be a maximal torus in $G$. If $M$ is a compact Hamiltonian $G$-space, and $0$
	is a regular value of the moment maps $\mu_{G}$ and $\mu_{T}$, then the $K$-theories of the symplectic quotients \(M \mod G\) and \(M \mod T\) are related by an isomorphism
\begin{equation}\label{eq:Ktheory-Martin}
K^*(M\mod G) \cong \frac{K^*(M \mod T)^W}{\Ann(e(E))},
\end{equation}
where $e(E)\in K^{*}(M\mod T)$ is the $K$-theoretic Euler class of the bundle $E$ given by (\ref{eq:bundle}).
\end{theorem}

In order to prove Theorems~\ref{theorem:martin} and \ref{theorem:Ktheory-Martin}, we must examine the Euler class $e(E)\in K^{*}(M\mod T)$ of the bundle $E = \mu_{T}^{-1}(0) \times_{T} \t^{0} \to M\mod T$ corresponding to its spin structure. We observe that while $E$ has a unique equivariant spin structure by Proposition~\ref{prop:topological-martin}, the sign of the Euler class $e(E)$ depends on a choice of orientation.
 However, as we are interested only in the annihilator of $e(E)$, this choice of sign does not affect Theorems~\ref{theorem:martin} or \ref{theorem:Ktheory-Martin}.
We first relate the $K$-theoretic Euler class $e(E)$ appearing in the statements of Theorems~\ref{theorem:martin} and~\ref{theorem:Ktheory-Martin} with the denominator of the Weyl character formula given in (\ref{eq:denominator}).

\begin{lemma}\label{lemma:euler}
	Given an $\Ad$-invariant inner product and a system of positive roots 
for $\g$, the Euler class $e(E)\in K^{*}(M\mod T)$ associated to the unique $T$-equivariant spin structure on $E$ is, up to sign, the image of the Weyl denominator $e_{T}(\g/\t)\in R(T)$ under the associated bundle map $R(T) \to K^{*}(M\mod T)$. 
\end{lemma}

\begin{proof}
Recalling Proposition~\ref{prop:topological-martin}, we can use the inner product on $\g$ to identify the annihilator $\t^{0}$ with the quotient $\g/\t$, and using the $\Ad$-invariance, the unique $T$-equivariant spin structure on $E$ is induced by the $T$-equivariant spin structure on $\g/\t$:
$$\Spin(E) \cong \mu_{T}^{-1}(0) \times_{T} \Spin(\g/\t).$$
It follows that the corresponding $K$-theoretic Euler class $e(E)$ is the image
of the $T$-equivariant Euler class $e_{T}(\g/\t)$ under the map
$R(T) \to K^{*}(M\mod T)$ taking virtual complex representations of $T$ to virtual complex vector bundles over $M\mod T$ associated to the principal $T$-bundle $\mu_{T}^{-1}(0)\to M\mod T$.

The Euler class associated to a spin structure is the virtual bundle given by the difference of the two corresponding complex half-spin bundles.
To construct the complex spin representation corresponding to a vector space $V$ we choose a polarization $V \otimes \C \cong W \oplus \overline{W}$, and the spin representation is then given by
$$\bbS_{V} := \Lambda^{*}_{\C}(\overline{W}) \otimes {\det}^{-1/2}(\overline{W}).$$
Here, the inverse square root of the determinant, i.e., the top exterior power, is a correction factor which renders this spin representation independent of the choice of polarization.
The corresponding $K$-theoretic
Euler class is the difference
$$e(V) = \Bigl( \left[\Lambda^{\text{even}}_{\C}(\overline{W})\right] - \
\left[\Lambda^{\text{odd}}_{\C}(\overline{W})\right] \Bigr)\otimes \left[{\det}^{-1/2}(\overline{W})\right].$$
Note that while the spin representation $\bbS_{V}$ is independent of the polarization, the sign of the Euler class $e(V)$, corresponding to the orientation of the vector space $V$, does depend on the choice of polarization.

Recall that the complement of the Cartan subalgebra decomposes as a direct sum of one-dimensional complex root spaces,
$$(\g/\t)\otimes \C \cong {\bigoplus}_{\alpha} \g_{\alpha}.$$
A $T$-invariant polarization of $\g/\t \otimes \C$ then corresponds to a system of positive roots for $\g$, with $W = \bigoplus_{\alpha>0}\g_{\alpha}$. Since the Euler class is multiplicative, we obtain
\begin{equation}\label{eq:e-gmodt}
e_{T}(\g/\t) = \prod_{\alpha>0}e_{T}(\g_{\alpha})
= \prod_{\alpha>0} \left(1 - e^{-\alpha}\right) e^{\alpha/2} 
= \prod_{\alpha>0} \left( e^{\alpha/2} - e^{-\alpha/2} \right)
\in R(T),\end{equation}
which is precisely the Weyl denominator (\ref{eq:denominator}).
\end{proof}

We note that the choice of positive root system in Lemma~\ref{lemma:euler} affects only the sign of the Euler class $e(E)$. We use this positive root system to determine the polarization of $E \otimes \C$ into isotropic subbundles, which in turn determines the orientation of $E$.

The main step in the proofs of Theorems~\ref{theorem:martin} and \ref{theorem:Ktheory-Martin} is the following lemma, which we prove using the results we established in Sections~\ref{sec:weyl} and \ref{sec:pushforward} above.

\begin{lemma}\label{lemma:main}
	For any Weyl-invariant element $y \in K^{*}(M\mod T)^{W}$, we have
	$$i^{*} y = 0 \Longleftrightarrow y \kcup e(E) = 0,$$
	or in other words
	$$\ker i^{*} \cap K^{*}(M\mod T)^{W} = \Ann\bigl(e(E)\bigr) \cap K^{*}(M\mod T)^{W}.$$
\end{lemma}

\begin{proof}
	Suppose that an element $y\in K^{*}(M\mod T)$ satisfies $i^{*}y = 0$. By Proposition~\ref{prop:topological-martin}, the map $i : \mu_{G}^{-1}(0)/T \to M\mod T$ is the inclusion of the zero set of a transverse section of the spin bundle $E\to M\mod T$. So, applying Lemma~\ref{lemma:transverse}, we find that
	$$0 = i_{!} (i^{*}y) = y \kcup e(E).$$
	
Conversely, suppose that $y \kcup e(E) = 0$ for a Weyl-invariant element $y \in K^{*}(M\mod T)^{W}$. Pulling back to $\mu_{G}^{-1}(0)/T$ via the inclusion $i$, we find that
$$0 = i^{*}y \kcup i^{*}e(E) \in K^{*}\bigl(\mu_{G}^{-1}(0)/T \bigr)
\cong K_{T}^{*}\bigl( \mu_{G}^{-1}(0)\bigr).$$
By Lemma~\ref{lemma:euler}, the class $e(E)\in K^{*}(M\mod T)$ is the image of the $T$-equivariant Euler class $e_{T}(\g/\t)$ under the map $R(T) \to K^{*}(M\mod T)$. Therefore, its pullback with respect to $i$ is likewise the image of $e_{T}(\g/\t)$ via the map
$$ e_{T}(\g/\t) \in R(T) \longmapsto i^{*}e(E) \in K^{*}_{T}\left(\mu_{G}^{-1}(0)\right). $$
By Corollary~\ref{cor:brion}, the element $i^{*}e(E)$ generates the
Weyl anti-invariants $K^{*}_{T}(\mu_{G}^{-1}(0))^{-W}$ as a free module over the Weyl-invariants $K_{T}^{*}(\mu_{G}^{-1}(0))^{W}$. In particular, since $i^{*}y$ is Weyl-invariant, we see that $i^{*}y$ must vanish.
\end{proof}

\begin{proof}[Proof of Theorem~\ref{theorem:martin}]
	Recall that in (\ref{eq:kernels}) we established
	$$\ker \kappa_{G} \cong \ker (i^{*} \circ \kappa_{T})|_{\im r^{G}_{T}}$$
	by studying the diagram (\ref{eq:Kirwan-GtoT-diagram}). The isomorphism
	between the left side and right sides of this equation is given by the
	restriction map $r^{G}_{T} : K_{G}^{*}(M) \to K_{T}^{*}(M)$.
	By Part (1) of Corollary~\ref{cor:brion}, we see that the image of
	the restriction map is precisely the Weyl-invariant elements,
	$$\im r^{G}_{T} = K_{T}^{*}(M)^{W},$$
	and so we have
	$$\ker\kappa_{G} \cong \left\{ x\in K_{T}^{*}(M)^{W} :
	i^{*} \kappa_{T}(x) = 0 \right\}.$$
	Since the Kirwan map $\kappa_{T}$ given by (\ref{eq:Kirwan-nonabelian}) is
	essentially the restriction map to $\mu_{T}^{-1}(0)$, we see that it commutes
	with the Weyl group action, and thus
	$$\kappa_{T} \bigl( K_{T}^{*}(M) \bigr) \subset K_{T}^{*}\left(\mu_{T}^{-1}(0)\right)^{W} \cong K^{*}(M\mod T)^{W}.$$
	So, in order to compute the kernel of $i^{*} \circ \kappa_{T}$ restricted to the Weyl-invariants $K_{T}^{*}(M)^{W}$, we must compute the kernel of $i^{*}$ restricted to the Weyl-invariants $K^{*}(M\mod T)^{W}$. However, by Lemma~\ref{lemma:main}, this kernel is precisely the annihilator of the Euler class $e(T)$, restricted to the Weyl-invariants $K^{*}(M\mod T)^{W}$, giving us
	$$\ker\kappa_{G} \cong \left\{ x\in K_{T}^{*}(M)^{W} : \kappa_{T}(x) \kcup e(E) = 0 \right \},$$
	which is our desired result.
\end{proof}

\begin{proof}[Proof of Theorem~\ref{theorem:Ktheory-Martin}]
	Referring back to the Weyl-invariant commutative diagram (\ref{eq:Kirwan-GtoT-diagram-Weyl}), we see that the composition
	 $$K^{*}(M\mod T)^{W} \xrightarrow{i^{*}} K^{*}\bigl( \mu_{G}^{-1}(0)/T\bigr)^{W}
	 \cong K^{*}(M\mod G)$$
is surjective, since $\kappa_G$ is surjective \cite{HL05}. We therefore have an isomorphism $K^{*}(M\mod G) \cong K^{*}(M\mod T)^{W} / \ker i^{*}$. Our result then follows directly from Lemma~\ref{lemma:main}, as the kernel of $i^{*}$ restricted to the Weyl-invariant elements is precisely the annihilator of $e(E)$.
\end{proof}

\section{Equal rank subgroups}

In this section, we consider generalizations of our main theorems where we replace the maximal torus $T$ of $G$ with a maximal rank subgroup $H$ of $G$. In other words, we consider subgroups $H\subset G$ for which there exists a common maximal torus $T$ such that $T \subset H \subset G$. For example, in the case $G = \U(n)$, the maximal rank subgroups are of the form $H = \U(n_{1}) \times \U(n_{2}) \times \cdots \times \U(n_{k})$ with $\sum n_{i} = n$. Another family of examples is $\SO(2n)$ inside $\SO(2n+1)$. The equal rank inclusion of $B4 = \Spin(9)$ inside the exceptional Lie group $F4$ has recently inspired some interesting representation theory \cite{GKRS98}.

 Let $W_{G}$ and $W_{H}$ be the corresponding Weyl groups. If in addition $W_{H}$ is a normal subgroup of $W_{G}$, we can then consider the relative Weyl group $W_{G,H} := W_{G} / W_{H}$, which does not depend on the choice of common maximal torus.
 In this situation, given a $G$-space $X$, the action of the Weyl group $W_{G}$ on $K^{*}_{T}(X)$ descends to an action of the relative Weyl group $W_{G,H}$ on $K^{*}_{H}(X) \cong K^{*}_{T}(X)^{W_{H}}$, in light of Corollary~\ref{cor:brion}. In particular, when $X$ is a point, the action of the Weyl group $W_{G}$ on $R(T)$ descends to an action  of the relative Weyl group $W_{G,H}$ on the representation ring $R(H) \cong R(T)^{W_{H}}$. We can now generalize our results of Section~\ref{sec:weyl} as follows:

\begin{lemma}
Let $G$ be a compact connected Lie group with no torsion in its fundamental group $\pi_{1}(G)$, and let $X$ be a compact $G$-space. If $H$ is a subgroup of maximal rank in $G$, then
\begin{equation}\label{eq:g-to-h}
K^{*}_{H}(X) \cong K^{*}_{G}(X) \otimes_{R(G)} R(H).
\end{equation}
Furthermore, if the pair $(G,H)$ admits a relative Weyl group $W_{G,H}$, then
\begin{enumerate}
\item
The restriction map $r^{G}_{H} : K_{G}^{*}(X) \to K_{H}^{*}(X)$ gives an isomorphism
$$r^{G}_{H} : K^{*}_{G}(X) \xrightarrow{\cong} K_{H}^{*}(X)^{W_{G,H}}.$$
\item
The $W_{G,H}$-anti-invariant elements $K_{H}^{*}(X)^{-W_{G,H}}$ are a free module over $K_{G}^{*}(X)$ via the restriction map $r^{G}_{H}$, with a single generator given by the image of the Euler class $e_{H}(\g/\h) \in R(H)^{-W_{G,H}}$, under the natural map $R(H)\to K_{H}^{*}(X)$.
\end{enumerate}
\end{lemma}

\begin{proof}
We can establish the identity (\ref{eq:g-to-h}) using the same argument as for Lemma~\ref{lemma:mcleod} above, considering here the K\"unneth spectral sequence for the $G$-equivariant $K$-theory of the \mbox{product} $X \times G/H$.  By our proof of Corollary~\ref{cor:brion}, we then need only verify the statements (1) and (2) when $X$ is a point, i.e., for the representation rings. For statement (1), we have
$$R(G) \cong R(T)^{W_{G}} \cong \left( R(T)^{W_{H}} \right)^{W_{G,H}}
\cong R(H)^{W_{G,H}}.$$

To establish statement (2), we first note that the relative Weyl group $W_{G,H}$ acts simply transitively on the Weyl chambers of $G$ contained inside a fixed Weyl chamber of $H$. The $W_{G,H}$-anti-invariant elements of $R(H)$ are therefore generated by principal alternating elements of the form
$$A([U]) := \sum_{w\in W_{G,H}} \sign(w) \,[U],$$ where $U$ is an irreducible representation of $H$. A choice of positive root system for $G$ determines a  positive root system for $H$. Letting $U_{\lambda}$ denote the irreducible representation of $H$ with highest weight $\lambda$, we see that $w[U_{\mu}] = [U_{\hat{w}\cdot \mu}]$, where $\hat{w}\in W_{G}$ is the lift of $w\in W_{G,H}$ which maps the positive Weyl chamber for $G$ into the positive Weyl chamber for $H$.
When considering the principal alternating elements generating $R(H)^{-W_{G,H}}$, we can restrict our attention to the additive basis of elements of the form
$$A(\lambda) := A\left( \left[ U_{\lambda + \rho_{G} - \rho_{H}} \right]\right),$$
where $\lambda$ is a dominant weight for $G$, and $\rho_{G}$ and $\rho_{H}$ are half the sums of the positive roots, or equivalently the sums of the basic weights, for $G$ and $H$, respectively. Adding $\rho_{G}$ shifts dominant weights for $G$ into the interior of the positive Weyl chamber for $G$, while subtracting $\rho_{H}$ shifts to include weights which lie on the border of the positive Weyl chamber for $H$.

We now recall the following generalization of the Weyl character formula due to Gross, Kostant, Ramond, and Sternberg \cite{GKRS98}:
\begin{equation}\label{eq:GKRS}
V_{\lambda} \otimes \bbS_{\g/\h}^{+} - V_{\lambda} \otimes \bbS_{\g/\h}^{-} = \sum_{w\in W_{G}/W_{H}}\sign(w)\,U_{\hat{w} (\lambda + \rho_{G})-\rho_{H}},
\end{equation}
where $V_{\lambda}$ is the irreducible representation of $\g$ with highest weight $\lambda$, and the identity is viewed in terms of virtual representations in $R(H)$. The difference of the two half-spin representations is precisely the $H$-equivariant $K$-theory Euler class of the representation $\g/\h$ of $H$:
\begin{equation}\label{eq:GKRS-denominator}
e_{H}(\g/\h) = \bigl[\bbS_{\g/\h}^{+}\bigr] - \bigl[\bbS_{\g/\h}^{-}\bigr] = \sum_{w\in W_{G}/W_{H}} \sign(w) \, [U_{\hat{w}\cdot \rho_{G} - \rho_{H}}] \,\in\, R(H).
\end{equation}
We note that (\ref{eq:GKRS}) and (\ref{eq:GKRS-denominator}) hold even when $W_{H}$ is not a normal subgroup of $W_{G}$. However, if $W_{H}$ is normal, these equations simplify slightly, since $\hat{w}\cdot \rho_{H} = \rho_{H}$ in this case, and thus
\begin{equation}\label{eq:rho-shift}
\hat{w}(\lambda + \rho_{G}) - \rho_{H} = \hat{w}(\lambda + \rho_{G} - \rho_{H}).\end{equation}
We can therefore rewrite (\ref{eq:GKRS}) in the form
$$A(\lambda) = r^{G}_{H} [V_{\lambda}] \cdot e_{H}(\g/\h),$$
and it follows that $R(H)^{-W_{G,H}}$ is a free module over $R(G)$, generated by $e_{H}(\g/\h)$.
\end{proof}

\begin{remark}We note that it is vital to this proof that $W_{H}$ be a normal subgroup of $W_{G}$. It not only allows us to consider actions of the relative Weyl group $W_{G,H} = W_{G} / W_{H}$, but also implies that the action of $\hat{w}$ fixes $\rho_{H}$, giving us (\ref{eq:rho-shift}) and allowing us to work with the relative $\rho$-shift $\rho_{G}-\rho_{H}$. This property holds, for example, if $G = G_{1} \times G_{2}$ and $H = G_{1} \times T_{2}$, where $G_{1}$ and $G_{2}$ are compact connected Lie groups, and $T_{2}$ is a maximal torus in $G_{2}$. In this case, we have $W_{G}\cong W_{G_{1}}\times W_{G_{2}}$, $W_{H} \cong W_{G_{1}}$, and $W_{G,H} \cong W_{G_{2}}$.
On the other hand, for many cases such as $G = \U(3)$ and $H = \U(2) \times \U(1)$,
we do not have a relative Weyl group, and this result does not apply.
\end{remark}

The results of Section~\ref{sec:pushforward} do not even mention the maximal torus $T$, and all that is required to apply them is an analogue of Martin's topological Proposition~\ref{prop:topological-martin}. However, the proof of Proposition~\ref{prop:topological-martin} does not require that $T$ be a maximal torus, or even abelian! So, we can replace $T$ in Proposition~\ref{prop:topological-martin} with any closed subgroup $H$ of $G$, maximal rank or otherwise.

Putting together all of these ingredients, the arguments of Section~\ref{sec:Ktheory-MmodG-MmodT} generalize \textit{mutatis mutandis}
to give the following theorem. Let $\kappa_G$ denote the Kirwan map as defined in~\eqref{eq:Kirwan-nonabelian}, and $\kappa_H$ the Kirwan map for the symplectic quotient $M \mod H$. Here, we replace the bundle $E$ given by (\ref{eq:bundle}) with the bundle
\begin{equation}\label{eq:bundle-h}
E_{H} := \mu_{H}^{-1}(0) \times_{H}\h^{0} \to M\mod H
\end{equation}
where $\h^{0}$ is the annihilator of $\h$ in $\g^{*}$.

\begin{theorem}
	Let $G$ be a compact connected Lie group with no torsion in its fundamental group $\pi_{1}(G)$, and let $H$ be a subgroup of maximal rank in $G$ which admits a relative Weyl group $W_{G,H}$. If $M$ is a compact Hamiltonian $G$-space, and $0$
	is a regular value of the moment maps $\mu_{G}$ and $\mu_{H}$, then the kernels of the Kirwan maps $\kappa_{G}$ and $\kappa_{H}$ given by (\ref{eq:Kirwan-nonabelian}) are related by
	$$\ker \kappa_{G}
	\cong \left\{ x\in K^{*}_{H}(M)^{W_{G,H}} : \kappa_{H}(x)\kcup e(E_{H}) = 0 \right\},$$
where $e(E_{H})\in K^{*}(M\mod H)$ is the $K$-theoretic Euler class of the bundle $E_{H}$ given by (\ref{eq:bundle-h}).
\end{theorem}

Finally, using our Kirwan surjectivity result from \cite{HL05}, we obtain our second stronger version of this theorem, computing the $K$-theory of the symplectic quotient rather than the kernel of the Kirwan map.

\begin{theorem}
	Let $G$ be a compact connected Lie group with no torsion in its fundamental group $\pi_{1}(G)$, and let $H$ be a subgroup of maximal rank in $G$ which admits a relative Weyl group $W_{G,H}$. If $M$ is a compact Hamiltonian $G$-space, and $0$
	is a regular value of the moment maps $\mu_{G}$ and $\mu_{H}$, then the $K$-theories of the symplectic quotients $M\mod G$ and $M \mod H$ are related by an isomorphism
	$$K^{*}(M\mod G) \cong \frac{K^{*}(M\mod H)^{W_{G,H}}}{\Ann\bigl(e(E_{H})\bigr)},$$
where $e(E_{H})\in K^{*}(M\mod H)$ is the $K$-theoretic Euler class of the bundle $E_{H}$ given by (\ref{eq:bundle-h}).
\end{theorem}

In both of these theorems, the bundle $E_{H}$ is associated to the
representation $\h^{0} \cong \g/\h$ of $H$, and the Euler class $e(E_{H})$
is induced by the $H$-equivariant Euler class $e_{H}(\g/\h) \in R(H)$
given in (\ref{eq:GKRS-denominator}) via the map
$$R(H) \to K_{H}^{*}\left(\mu_{H}^{-1}(0)\right) \cong K^{*}(M\mod
H).$$ This Euler class can also be viewed as the $G$-equivariant Euler
class $e_{G}(G/H)$ via the isomorphism $K_{G}^{*}(G/H)\cong R(H)$, as
well as the ``denominator'' in the maximal rank generalization
(\ref{eq:GKRS}) of the Weyl character formula. For further discussion
of this Euler class, see Bott's classic papers \cite{Bott65, Bott88}, or the
second author's recent work \cite{Lan-twisted}. For further discussion
of the formula (\ref{eq:GKRS}) and the ``Euler number multiplets''
which appear on its right hand side, see Kostant's paper \cite{Kos99}
or the second author's extension of these results to loop groups
\cite{Lan01}.

\section{Other cohomology theories}

In this section, we discuss versions of our main results, Theorems~\ref{theorem:martin} and~\ref{theorem:Ktheory-Martin}, in different cohomology theories. We begin with rational cohomology, comparing our treatment with that of Martin in \cite{Mar00}. We then discuss integral cohomology and, finally, other 
cohomology theories such as complex cobordism. 



Since the Chern character gives a ring isomorphism from rational $K$-theory to rational cohomology, our $K$-theoretic results imply the corresponding statements in rational cohomology. We note that the Chern character of a $K$-theory Euler class is not equal to the corresponding cohomology Euler class; rather, their quotient is the Todd class. However, the Todd class is invertible, and thus the two have the same annihilator ideals.

We present two different methods for proving the rational cohomology analogues of our Theorems~\ref{theorem:martin} and~\ref{theorem:Ktheory-Martin} directly.
The first method is to use direct rational-cohomological analogues of the two vital ingredients in our proofs in Section~\ref{sec:Ktheory-MmodG-MmodT}, namely Corollary~\ref{cor:brion} and Lemma~\ref{lemma:transverse}. The analogous results both hold for rational cohomology, as we discuss below.  The second method considers an additional pushforward map corresponding to the fibration \(\pi: \mu_G^{-1}(0)/T \to M \mod G.\) This latter method is more in the spirit of Martin's original argument, but still differs from it in some respects. We discuss the differences in detail below. 

 We begin with the direct rational cohomology analogue of our proofs. As an analogue to part (1) of Corollary~\ref{cor:brion}, 
 it is well known that if $X$ is a
$G$-space, then the restriction map gives an isomorphism
$$r^{G}_{T} : H_{G}^{*}(X;\Q) \xrightarrow{\cong} H^{*}_{T}(X;\Q)^{W}$$
between the $G$-equivariant cohomology and the Weyl-invariant part of the $T$-equivariant cohomology.
Moreover, in \cite{Bri91}, Brion argues that the Weyl anti-invariant part $H^{*}_{T}(X;\Q)^{-W}$ is a free module of rank $1$ over $H^{*}_{G}(X;\Q)$, and  Brion further notes that, via the decomposition
$$H^{*}_{T}(X;\Q) \cong H^{*}_{G}(X;\Q) \otimes_{\Q} H^{*}(G/T;\Q),$$
the single generator is induced by any non-vanishing element of the top cohomology group $$H^{*}(G/T;\Q)^{-W} \cong H^{\text{top}}(G/T;\Q).$$
In particular, the cohomology Euler class $e^{H}(G/T) \in H^{*}(G/T;\Q)$ is such an element,
and so $1 \otimes e^{H}(G/T)$ generates $H^{*}_{T}(X;\Q)^{-W}$ as a module over $H^*_G(X;\Q)$. However, this Euler class is induced by the equivariant Euler class $e^{H}_{T}(\g/\t)\in H^{*}_{T}(\pt;\Q)$ via the natural homomorphism
$$H^{*}_{T}(\pt;\Q) \to H_{T}^{*}(G;\Q) \cong H^{*}(G/T;\Q),$$
and thus the image of $e^{H}_{T}(\g/\t)$ generates $H^{*}_{T}(X;\Q)^{-W}$ as a free module over $H^{*}_{G}(X;\Q)$. This is the rational cohomology analogue of part (2) of our Corollary~\ref{cor:brion}. 

The analogue of Lemma~\ref{lemma:transverse} in de~Rham cohomology is the
statement that if $E \to Y$ is an oriented vector bundle and a section
$s : Y \to E$ is transverse to the zero section, then the zero set of
$s$ is Poincar\'{e} dual to the cohomology Euler class $e^{H}(E)\in H^{\rank
E}(X)$ (see \cite{BT82}). In the special case where $E$ is the
tangent bundle of an oriented manifold, this becomes the familiar fact
that the Euler characteristic is the number of zeros of a generic
vector field. We note that our proof of Lemma~\ref{lemma:transverse}
is not specific to $K$-theory, and in fact holds for 
rational or de~Rham cohomology when $E$ is an oriented vector bundle.

Putting together these prerequisites, the argument given in
Section~\ref{sec:Ktheory-MmodG-MmodT} can now be used to give a proof
of the rational cohomology versions of Theorems~\ref{theorem:martin}
and \ref{theorem:Ktheory-Martin}. In contrast, Martin's proof in
\cite{Mar00} does not consider anti-invariant elements, bypassing
Brion's lemma entirely. Instead, Martin considers the pushforward map
$$\pi_{!} : H^{*}\bigl(\mu_{G}^{-1}(0)/T;\Q\bigr) \to H^{*}(M\mod
G;\Q)$$ corresponding to integration along the $G/T$ fibers of the map
$\pi : \mu_{G}^{-1}(0)/T \to M\mod G$ given in
(\ref{eq:topological-Martin}). In a statement analogous to that of our
Lemma~\ref{lemma:transverse}, Martin argues that if $\vertical$ is the
bundle of vertical tangent vectors over $\mu_{G}^{-1}(0)/T$, then
\begin{equation}\label{eq:martin-step2}
\pi_{!} \bigl( e^{H}(\vertical) \hcup \pi^{*}x \bigl) = \pi_{!}\bigl(e^{H}(\vertical)\bigr) \hcup x = |W|\,x
\end{equation}
for all $x \in H^{*}(M\mod G;\Q)$. 

Using this additional pushforward map $\pi_{!}$ in~\eqref{eq:martin-step2}, we now give our second proof of the rational cohomology analogue of our fundamental Lemma~\ref{lemma:main}.
Let $e^{H}(E) \in H^*(M \mod T;\Q)$ be the cohomology Euler class of the bundle \(E \to M \mod T\) given by \eqref{eq:bundle}. This is the image of the equivariant Euler class $e_T(\g/\t)$ under the map
\[H^{*}_T(\pt;\Q) \to H^{*}_T\left(\mu_T^{-1}(0);\Q\right) \cong H^{*}(M \mod T;\Q).\] 

\begin{lemma}\label{lemma:martin}
For any Weyl-invariant element $y\in H^{*}(M\mod T;\Q)^{W}$, we have
$$i^{*}y = 0 \Longleftrightarrow y \hcup e^{H}(E) = 0,$$
where $e^{H}(E)$ is induced by the equivariant cohomology Euler class $e^{H}_{T}(\g/\t)\in H^{*}_{T}(\pt;\Q)$.
\end{lemma}

\begin{proof}
Suppose that an element $y\in H^{*}(M\mod T;\Q)$ satisfies $i^{*}y = 0$. By the rational cohomology analogue of Lemma~\ref{lemma:transverse}, we have
$$0 = i_{!}(i^{*}y) = y \hcup e^{H}(E).$$
Conversely, suppose that $y \hcup e^{H}(E) = 0$ for a Weyl-invariant element $y\in H^{*}(M\mod T;\Q)^{W}$. Pulling back to $\mu_{G}^{-1}(0)/T$ via the inclusion $i$, we have
\begin{equation}\label{eq:martin-step}
0 = i^{*}y \hcup i^{*}e^{H}(E) \in H^{*}\bigl(\mu_{G}^{-1}(0)/T;\Q\bigr).
\end{equation}
Since $y$ is Weyl-invariant, so is its pullback $i^{*}y$. However, since the map $$\pi^{*} : H^{*}(M\mod G;\Q) \to H^{*}\bigl(\mu_{G}^{-1}(0)/T;\Q\bigr)^{W}$$ is isomorphic to the restriction
$$r^{G}_{T} : H^{*}_{G}\bigl(\mu_{G}^{-1}(0);\Q\bigr) \xrightarrow{\cong} H^{*}_{T}\bigl(\mu_{G}^{-1}(0);\Q\bigr)^{W},$$
we see that $i^{*}y = \pi^{*}x$ for some $x\in H^{*}(M\mod G;\Q)$. Applying the pushforward $\pi_{!}$ to both sides of (\ref{eq:martin-step}) and using (\ref{eq:martin-step2}), we obtain
$$0 = \pi_{!} \bigl(i^{*}y \hcup i^{*}e^{H}(E) \bigr)
= \pi_{!}\bigl( \pi^{*}x \hcup e^{H}(\vertical) \bigr)
= |W|\, x,$$
where we note that the bundles $i^{*}E$ and $\vertical$ over $\mu_{T}^{-1}(0)/T$ are induced by the isomorphic representations $\t^{0}$ and $\g/\t$ of $T$. Since $|W|$ is strictly positive, we must have $x =0$, and therefore $i^{*}y = \pi^{*}x = 0$.
\end{proof}

The rational cohomology analogues of our main theorems follow immediately from the above lemma, following the arguments in Section~\ref{sec:Ktheory-MmodG-MmodT}. 
We note that while this lemma is attributed to Martin in \cite{TW03},
it does not actually appear in his paper \cite{Mar00}. In
particular, all of the results in \cite{Mar00} are given in terms of
the square of the cohomology Euler class $e^{H}(E)^{2}$, whereas ours, including the lemma above, are in terms of $e(E)$. 
Nevertheless, our
proof above of Lemma~\ref{lemma:martin} captures the spirit of
Martin's arguments in \cite{Mar00}, since we examine, as he does, the behavior
of the cohomology Euler class $e^{H}(E)$ with respect to both of the pushforwards $i_!$ and $\pi_!$.
 However, we have removed from Martin's arguments any reliance
on Poincar\'{e} duality pairings and averaging over the Weyl group.
Similarly, in \cite{JMW02}, Jeffrey, Mare, and Woolf prove the
equivalence of the two results, one stated in terms of the single cohomology Euler class $e^{H}(E)$ and
the other using its square $e^{H}(E)^2$. They follow Martin's general setup, and thus use a summation over the Weyl group and Poincar\'e duality, in addition to Brion's lemma. However, as we
showed in Section~\ref{sec:Ktheory-MmodG-MmodT}, once we consider
Brion's lemma, these results follow {\em directly} without use of integral pairings or Poincar\'e duality.

If we were to consider integral cohomology rather than rational
cohomology, we encounter several problems.  In order for some of
the results analogous to those in Section~\ref{sec:weyl} to hold for $H^*(-;\Z)$, we must require that the cohomology $H^{*}(G;\Z)$ of the group $G$ be torsion-free. This is a much more restrictive condition than the $K$-theoretic requirement that $\pi_{1}(G)$ be torsion-free, and in particular Borel showed in
\cite{Bor54} that the Lie groups $\Spin(n)$ for $n \geq 7$, as well as the exceptional Lie groups $G_{2}$ and $F_{4}$, all have torsion in their integral cohomology. 
More importantly, in the integral case, part (2) of
Corollary~\ref{cor:brion} does not hold as written. It is true that the Weyl anti-invariants are a free module of rank 1 over
the Weyl-invariants, but the equivariant cohomology Euler class $e^{H}_{T}(\g/\t)$ is not a
generator! Rather, the generator is $1/|W|\,e^{H}_{T}(\g/\t)$, which we
see by noting that $e^{H}(G/T)$ is $\chi(G/T) = |W|$ times the generator
of the top cohomology of $G/T$. As a result, we cannot repeat our proof of
Lemma~\ref{lemma:main}, since in the presence of torsion the
element $|W|$ may become a zero divisor. Similarly, we cannot use the argument in our proof of  Lemma~\ref{lemma:martin} above, since that also requires division by $|W|$. Finally, to prove an integral cohomology version of
Theorem~\ref{theorem:Ktheory-Martin}, we would require an integral cohomology version of Kirwan surjectivity,
which does not necessarily hold because of the presence
of torsion (see \cite{TW03}).

The difficulties in integral cohomology may be eliminated by inverting those primes which appear in the torsion of $H^{*}(G;\Z)$, which divide the order $|W|$ of the Weyl group, or which interfere with the Kirwan surjectivity theorem. However, none of this is necessary in $K$-theory.
We therefore conclude that (integral) $K$-theory, rather than integral
cohomology, is the natural setting in which to extend these rational
cohomology results. This agrees with our observation in \cite{HL05}
that passing from integral cohomology to $K$-theory eliminates
precisely enough torsion for these results to hold. Finally, we conjecture that
these results also have natural extensions to complex cobordism,
which determines complex $K$-theory much as complex $K$-theory
determines rational cohomology. For instance, in \cite{Hod75}, Hodgkin
derives his K\"unneth spectral sequence for any cohomology theory, and
we expect that results analogous to those in Section~\ref{sec:weyl}
hold for complex cobordism as well. Furthermore, the results of
Section~\ref{sec:pushforward} hold in any cohomology theory where the
bundles admit Thom isomorphisms and Euler classes. In particular, the
bundle $E$ associated to the representation $\g/\t$ admits a complex
structure, and so it gives rise to a corresponding Thom isomorphism
and Euler class in complex cobordism. We also believe it worthwhile to
explore the versions of these results in other variants of $K$-theory
and cobordism, such as $KO$-theory or spin cobordism, as well as other
cohomology theories such as elliptic cohomology.


\begin{thebibliography}{10}

\bibitem{Ati68}
M.~F. Atiyah.
\newblock Bott periodicity and the index of elliptic operators.
\newblock {\em Quart. J. Math. Oxford Ser. (2)}, 19:113--140, 1968.

\bibitem{Ati70}
M.~F. Atiyah.
\newblock Global theory of elliptic operators.
\newblock In {\em Proc. Internat. Conf. on Functional Analysis and Related
  Topics (Tokyo, 1969)}, pages 21--30. Univ. of Tokyo Press, Tokyo, 1970.

\bibitem{AB82}
M.~F. Atiyah and R.~Bott.
\newblock The {Y}ang-{M}ills equations over {R}iemann surfaces.
\newblock {\em Philos. Trans. Roy. Soc. London Ser. A}, 308(1505):523--615,
  1983.

\bibitem{AS68I}
M.~F. Atiyah and I.~M. Singer.
\newblock The index of elliptic operators. {I}.
\newblock {\em Ann. of Math. (2)}, 87:484--530, 1968.

\bibitem{Bor54}
A.~Borel.
\newblock Sur l'homologie et la cohomologie des groupes de {L}ie compacts
  connexes.
\newblock {\em Amer. J. Math.}, 76:273--342, 1954.

\bibitem{BS58}
A.~Borel and J.-P. Serre.
\newblock Le th\'eor\`eme de {R}iemann-{R}och (d'apr\`es {A}. {G}rothendieck).
\newblock {\em Bull. Soc. Math. France}, 86:97--136, 1958.

\bibitem{Bott65}
R.~Bott.
\newblock The index theorem for homogeneous differential operators.
\newblock In {\em Differential and Combinatorial Topology (A Symposium in Honor
  of Marston Morse)}, pages 167--186. Princeton University Press, Princeton,
  N.J., 1965.

\bibitem{Bott88}
R.~Bott.
\newblock On induced representations.
\newblock In {\em The mathematical heritage of Hermann Weyl (Durham, NC,
  1987)}, volume~48 of {\em Proc. Sympos. Pure Math.}, pages 1--13. Amer. Math.
  Soc., Providence, RI, 1988.

\bibitem{BT82}
R.~Bott and L.~W. Tu.
\newblock {\em Differential forms in algebraic topology}, volume~82 of {\em
  Graduate Texts in Mathematics}.
\newblock Springer-Verlag, New York, 1982.

\bibitem{Bri91}
M.~Brion.
\newblock Cohomologie \'equivariante des points semi-stables.
\newblock {\em J. Reine Angew. Math.}, 421:125--140, 1991.

\bibitem{Gol02}
R.~F. Goldin.
\newblock An effective algorithm for the cohomology ring of symplectic
  reductions.
\newblock {\em Geom. Funct. Anal.}, 12(3):567--583, 2002, math.SG/0110022.

\bibitem{GKM}
M.~Goresky, R.~Kottwitz, and R.~MacPherson.
\newblock Equivariant cohomology, {K}oszul duality, and the localization
  theorem.
\newblock {\em Invent. Math.}, 131:25--83, 1998.

\bibitem{GKRS98}
B.~Gross, B.~Kostant, P.~Ramond, and S.~Sternberg.
\newblock The {W}eyl character formula, the half-spin representations, and
  equal rank subgroups.
\newblock {\em Proc. Natl. Acad. Sci. USA}, 95(15):8441--8442 (electronic),
  1998, math.RT/9808133.

\bibitem{HL-kernel}
M.~Harada and G.~D. Landweber.
\newblock The ${K}$-theory of symplectic quotients.
\newblock In preparation.

\bibitem{HL05}
M.~Harada and G.~D. Landweber.
\newblock Surjectivity for {H}amiltonian {$G$}-spaces in {$K$}-theory.
\newblock {\em Trans. Amer. Math. Soc.}, to appear, math.SG/0503609.

\bibitem{Hat02}
A.~Hatcher.
\newblock {\em Algebraic topology}.
\newblock Cambridge University Press, Cambridge, 2002.

\bibitem{HP05}
T.~Hausel and N.~Proudfoot.
\newblock Abelianization for hyperk\"ahler quotients.
\newblock {\em Topology}, 44(1):231--248, 2005, math.SG/0310141.

\bibitem{Hod75}
L.~Hodgkin.
\newblock The equivariant {K}\"{u}nneth theorem in {$K$}-theory.
\newblock In {\em Topics in $K$-theory. Two independent contributions}, volume
  496 of {\em Lecture Notes in Math.}, pages 1--101. Springer, Berlin, 1975.

\bibitem{JK95}
L.~C. Jeffrey and F.~C. Kirwan.
\newblock Localization for nonabelian group actions.
\newblock {\em Topology}, 34:291--327, 1995.

\bibitem{JMW02}
L.~C. Jeffrey, A.-L. Mare, and J.~M. Woolf.
\newblock The {K}irwan map, equivariant {K}irwan maps, and their kernels.
\newblock {\em J. Reine Angew. Math}, to appear, math.SG/0211297.

\bibitem{Kir84}
F.~Kirwan.
\newblock {\em Cohomology of quotients in symplectic and algebraic geometry},
  volume~31 of {\em Mathematical Notes}.
\newblock Princeton University Press, Princeton, N.J., 1984.

\bibitem{Kir98}
F.~Kirwan.
\newblock Momentum maps and reduction in algebraic geometry.
\newblock {\em Differential Geom. Appl.}, 9(1-2):135--171, 1998.

\bibitem{Kos99}
B.~Kostant.
\newblock A cubic {D}irac operator and the emergence of {E}uler number
  multiplets of representations for equal rank subgroups.
\newblock {\em Duke Math. J.}, 100(3):447--501, 1999.

\bibitem{Lan01}
G.~D. Landweber.
\newblock Multiplets of representations and {K}ostant's {D}irac operator for
  equal rank loop groups.
\newblock {\em Duke Math. J.}, 110(1):121--160, 2001, math.RT/0005057.

\bibitem{Lan-twisted}
G.~D. Landweber.
\newblock Twisted representation rings and {D}irac induction.
\newblock {\em J. Pure Appl. Algebra}, to appear, math.RT/0403524.

\bibitem{Mar00}
S.~Martin.
\newblock Symplectic quotients by a nonabelian group and by its maximal torus,
  January 2000, math.SG/0001002.

\bibitem{McL79}
J.~McLeod.
\newblock The {K}\"{u}nneth formula in equivariant {$K$}-theory.
\newblock In {\em Algebraic topology, Waterloo, 1978 (Proc. Conf., Univ.
  Waterloo, Waterloo, Ont., 1978)}, volume 741 of {\em Lecture Notes in Math.},
  pages 316--333. Springer, Berlin, 1979.

\bibitem{Mer98}
A.~S. Merkur'ev.
\newblock Comparison of the equivariant and the ordinary {$K$}-theory of
  algebraic varieties.
\newblock {\em St. Petersburg Math. J.}, 9:815--850, 1998.

\bibitem{Pit72}
H.~V. Pittie.
\newblock Homogeneous vector bundles on homogeneous spaces.
\newblock {\em Topology}, 11:199--203, 1972.

\bibitem{Seg68}
G.~Segal.
\newblock Equivariant {$K$}-theory.
\newblock {\em Inst. Hautes \'Etudes Sci. Publ. Math.}, 34:129--151, 1968.

\bibitem{Sna72}
V.~P. Snaith.
\newblock On the {K}\"{u}nneth formula spectral sequence in equivariant
  {$K$}-theory.
\newblock {\em Proc. Cambridge Philos. Soc.}, 72:167--177, 1972.

\bibitem{TW03}
S.~Tolman and J.~Weitsman.
\newblock The cohomology rings of symplectic quotients.
\newblock {\em Comm. Anal. Geom.}, 11(4):751--773, 2003, math.DG/9807173.

\bibitem{VeVi03}
G.~Vezzosi and A.~Vistoli.
\newblock Higher algebraic {$K$}-theory for actions of diagonalizable groups.
\newblock {\em Invent. Math.}, 153(1):1--44, 2003, math.AG/0107174.

\end{thebibliography}

\def\cprime{$'$}

\end{document}